\documentclass[12pt]{amsart}
\usepackage{amssymb}

\input amssym.def

 2

\newcommand{\R}{{\bf R}}
\newcommand{\C }{{\bf C}}
\newcommand{\Z}{{\bf Z}}
\newcommand{\Q}{{\bf Q}}
\newcommand{\N}{{\bf N}}
\newcommand{\SL}{{\bf SL}}
\def\bfminus{\mathrel{\hbox{\bf -\kern-1pt-}}}
\def\ckomma{\, \raise2pt \hbox{,}}

\begin{document}

\title[Modular Forms]{The circle method and non lacunarity of modular functions}

\author{Sanoli Gun and Joseph Oesterl\'e}

\address[Sanoli Gun]
      {The Institute of Mathematical Sciences,
      C.I.T Campus, Taramani, Chennai 600113,
       India.}
\address[Joseph Oesterl\'e]
      {Institut de Math\'ematiques de Jussieu,
       175 rue du Chevaleret,  
       75013-Paris, 
       France.}
\email[Sanoli Gun]{sanoli@imsc.res.in}
\email[J. Oesterl\'e]{oesterle@math.jussieu.fr}

\renewcommand{\thefootnote}{}
\footnote{ \noindent\textbf{} \vskip0pt
\textbf{2010 Mathematics Subject Classification:} 
{Primary 11F20; Secondary 11F30, 11F37.}}

\maketitle

\begin{abstract}
Serre proved that any holomorphic cusp form of weight one for $\Gamma_1({\rm N})$
is lacunary while a holomorphic modular form for $\Gamma_1({\rm N})$  
of higher integer weight is lacunary if and only if it is a linear
combination of cusp forms of CM-type (see \cite{serre}, subsections 7.6 and 7.7). 
In this paper, we show that when a non-zero modular function of 
arbitrary real weight for any finite index subgroup of the modular 
group ${\SL}_2(\Z)$  is lacunary, it is necessarily holomorphic on 
the upper-half plane, finite at the cusps and has non-negative weight. 
\end{abstract}

\section*{Introduction}

Let $g$ be a holomorphic function on some half plane
$\goth{H}_r = \left\{z \in \C ~|~ {\rm Im}(z) > r \right\}$
with $r \in \R$. Assume that there exists some real number $h >0$
and some complex number $u$ such that $|u|=1$ and $g(z+h)=ug(z)$
for $z\in \goth{H}_r$. Then $g$ has a unique expansion of the form
\begin{eqnarray}\label{eq1}
g(z)= \sum_{\lambda \in \Lambda} a(\lambda) e^{2\pi i\lambda z}, 
\end{eqnarray}
where $\Lambda$ consists of all $\lambda \in \R$ such that
$e^{2\pi i\lambda h}= u$. We have $\Lambda = \lambda_0 + \Z h^{-1}$ for any 
$\lambda_0 \in \Lambda$.
The series \eqref{eq1} is called the {\em expansion of $g$
at $\infty$}. It converges normally on any horizontal strip 
$\left\{z \in \C ~|~ r' < {\rm Im}(z) < r'' \right\}$
with $r < r' < r''$. For $ \lambda \in \Lambda$, the coefficient $a(\lambda)$ is given by
\begin{eqnarray}\label{eq2}
a(\lambda) = h^{-1} \int_{[\tau, ~\tau+ h]}g(z)e^{-2\pi i\lambda z} dz
\end{eqnarray} 
where $\tau$ is any point of $\goth{H}_r$. We shall say that: \par 
-- $g$ {\em has finite order at $\infty$} if
$\left\{\lambda <0 ~|~ a(\lambda)\ne 0 \right\}$ is finite;\par
-- $g$ {\em is finite at $\infty$} if $a(\lambda)=0$ for all $\lambda<0$;
in that case \eqref{eq1} converges normally on any half plane 
$\goth{H}_{r'}$ with $r'> r$.\par
  
If $g$ has finite order at $\infty$, we shall say that its 
expansion at~$\infty$:\par
-- is {\em lacunary} if $\#\left\{\lambda \leqslant x ~|~ a(\lambda)\ne 0 \right\}= 
{\rm o}(x)$
when $x \to +\infty$;\par
-- is {\em strongly non-lacunary} if there exists a non-constant arithmetic progression 
$\Lambda' = \lambda_0 + m h^{-1}\N$ contained in $\Lambda$ 
such that $a(\lambda) \ne 0$ for all $\lambda \in \Lambda'$; in that case 
the expansion of $g$ at $\infty$ is not lacunary.\par

All these definitions depend only on $g$ and not on the choice of $h$
and $u$.

\bigskip

We now consider a finite index subgroup $\Gamma$ of $\SL_2(\Z)$ and a real
number~$k$. We define $z^k$ for $z \in \C^{*}$ to be $e^{k \log z}$,
where the logarithm of $z$ is chosen with imaginary part 
in $]-\pi, \pi]$. 

\bigskip
\noindent
{\bf Definition 1.}$-$
{\em We say that a meromorphic function $f$ on the upper half plane 
$\goth{H}= \left\{z \in \C ~|~ {\rm Im}(z) > 0 \right\}$ is 
{\em a modular function of weight $k$ for} $\Gamma$ if: \par
$a)$ for each $\gamma = \big({ a\ b \atop c\ d }\big) \in \Gamma$, $f$  
satisfies a functional equation
of the form $f(\gamma z)= j(\gamma, z)f(z)$, where the 
automorphy factor $j(\gamma,z)$ is a
holomorphic function on $\goth H$ such that 
$|j(\gamma, z)|= |cz+d|^k$ and where $\gamma z= \frac{az+b}{cz+d}${\em ;}\par
$b)$
for any $ \big({ a\ b \atop c\ d }\big) \in \SL_2(\Z)$, the function 
$z \mapsto (cz+d)^{-k}f(\frac{az+b}{cz+d})$ is holomorphic on some 
half plane ${\goth H}_r$
{\em (hence by $a)$ satisfies the conditions stated for $g$ at the 
beginning of this introduction)}
and it has finite order at $\infty$.}

\medskip
Note that when $f \ne 0$, condition $a)$ determines the 
automorphy factors $j(\gamma,z)$.
They satisfy the {\em cocycle relation} 
\begin{equation}\label{cocycle}
j(\gamma \gamma', z) = j(\gamma, \gamma' z)j(\gamma', z)
\end{equation}
for $\gamma, \gamma' \in \Gamma$ and $z \in \goth{H}$. Condition $b)$ is then
expressed by saying that {\em $f$ has finite order at the cusps}; if moreover 
$z \mapsto (cz+d)^{-k}f(\frac{az+b}{cz+d})$ is finite at $\infty$ for
any $ \big({ a\ b \atop c\ d }\big) \in \SL_2(\Z)$, then $f$ is said to be 
{\em finite at the cusps}. 

\bigskip
\noindent
{\bf Definition 2.}$-$
{\em A modular function is said to be {\rm lacunary 
(resp. strongly non-lacunary)} if its 
expansion at $\infty$ is lacunary {\rm (resp.} strongly non-lacunary{\rm)}.}

\bigskip
We can now state the main theorems of this paper. Let $k$ be a 
real number
and $\Gamma$  be a finite index subgroup of $\SL_2(\bf Z)$. Let 
$f $ be a non-zero modular function
of weight $k$ for $\Gamma$, with respect to some automorphy factors.

\bigskip
\noindent{\bf Theorem 1.}$-$
{\em If $f$ is not holomorphic on $\goth{H}$, $f$  is not lacunary.} 

\bigskip
From now on, we assume $f$ to be holomorphic on $\goth{H}$.

\bigskip
\noindent{\bf Theorem 2.}$-$
{\em If $k<0$, $f$ is strongly non-lacunary.}

\bigskip
\noindent{\bf Theorem 3.}$-$
{\em If $k \geqslant 0$ and $f$ is not finite at the cusps, $f$ is 
strongly non-lacunary.}

\bigskip
In the case $\Gamma=\SL_2(\Z)$, one can be more precise. There 
exists $u \in \C^{*}$ with $|u| = 1$ such that $f(z+1) = uf(z)$ and the expansion of $f$ at 
$\infty$ takes the form
$\sum_{\lambda \in \Lambda} a(\lambda) e^{2\pi i\lambda z}$, where 
$\Lambda = \left\{\lambda \in \C |  e^{2\pi i\lambda} = u\right\}$ . 
The conclusion that ``$f$ is strongly non-lacunary''
of theorem 2 and theorem 3 can then
be strengthened as follows: there exists $\lambda_0 \in \R$ such that
$a(\lambda) \ne 0$ for all $\lambda \geqslant \lambda_0$ in $\Lambda$. Moreover, 
such a $\lambda_0$ can be determined effectively.

\bigskip
Theorem 1 is not very deep and has just been stated for the sake of completeness. 
In fact any
meromorphic function $g$ on $\goth{H}$ satisfying the conditions stated at 
the beginning of the introduction, having finite order at $\infty$  and
whose expansion at $\infty$ is lacunary, is necessarily holomorphic on $\goth{H}$. 
This will be proved in section~5.

\medskip
Theorem 2 will be deduced from the Hardy-Ramanujan circle method (see \cite{HR}), 
as generalized by
Rademacher (see \cite{rad}, Chapters 14 and 15). We shall present this theory in 
section 1 in the special case where
$\Gamma= \SL_2(\Z)$: if $f$ is a modular function for $\SL_2(\Z)$ of some real 
weight $k<0$ for some
automorphy factors, Rademacher's theory yields convergent series expressing 
the coefficients of the expansion of $f$ at $\infty$ in terms of those
of its polar part (and in particular closed formulas for the partition function $p(n)$ 
when $f$ is the inverse of the Dedekind
eta function). A variant of this theory leading to the same results but using 
Poincar\'e
series will be presented in section 2.  

\medskip
Rademacher's theory can be extended to the case where $k \geqslant 0$ except
 that in this case
the formulas for the coefficients involve an error term. This will be explained in
section 3 and theorem 3 will follow in the case $\Gamma= \SL_2(\Z)$. 

\medskip
In section 4, we shall explain how to extend theorems 2 and 3 to any finite 
index subgroup $\Gamma$ of $\SL_2(\Z)$.

\medskip
The motivation of this paper arose from \cite{CGR} and \cite{GR}, where 
classifications of the lacunary
eta-products of some particular types were given. All those found had a 
non-negative weight and 
were finite at the cusps. Our goal was to find a theoretical explanation for
this observation. This explains the degree of generality (arbitrary real weight,
automorphy factors, finite index subgroup of $\SL_2(\Z)$) chosen to write this paper.

\section*{{\S 1.} The circle method for negative weight, \\
                following Rademacher}

\subsection*{1.0}{\bf Assumptions and notations} 

\smallskip
In this section, we consider a non-zero modular function $f$ for $\SL_2(\Z)$ 
of some real 
weight $k$,
{\em holomorphic on $\goth H$}. We denote by $j(\gamma,z)$ its automorphy factors. 
Hence we have 
\begin{eqnarray}\nonumber
f(\gamma z) = j(\gamma,z) f(z)  \quad \hbox{and} \quad |j(\gamma,z)| = |cz+d|^k  
\end{eqnarray}
for  $\gamma = \big({ a\ b \atop c\ d }\big)$ in $\SL_2(\Z)$  and $z$ in $\goth{H}$. 
Moreover, the functions $j(\gamma,z)$ satisfy the cocycle relation 
$j(\gamma \gamma', z) = j(\gamma, \gamma' z)j(\gamma', z)$.
The function $j(\big({ 1\ 1 \atop 0\ 1 }\big) , z)$ takes a constant value $u$.
We have $f(z+1) =uf(z)$ for $z \in \goth H$ and $|u|=1$. Let 
$$
f(z)=\sum_{\lambda \in \Lambda} a(\lambda) e^{2 \pi i\lambda z}
$$
be the expansion of $f$ at $\infty$,
where $\Lambda$ is the set of $\lambda \in \R$ such that $e^{2\pi i\lambda} =u$. 
By definition of
a modular function there are only finitely many $\lambda < 0$ in $\Lambda$ such that
$a(\lambda) \ne 0$.  We call 
$$
{\rm P}(z)= \sum_{\lambda \in \Lambda, \lambda<0} a(\lambda)e^{2\pi i \lambda z},
$$ 
the {\em polar part of $f$ at $\infty$}.

\subsection*{1.1}{\bf The circle method} 

\smallskip
For each rational number $x$ written in an irreducible form $\frac{a}{c}$, 
let ${\rm C}_x$ denote the circle in the upper half plane
of radius $1/2c^2$ tangent to the real line at the point $x$. 
These circles are called the {\em Ford circles}. We orient them in 
the clockwise direction and denote by ${\rm C}_x^{*}$ the circle
${\rm C}_x$ with the point $x$ removed.

\bigskip
\noindent
{\bf Proposition 1.} $-$
{\em Assumptions and notations are those of subsection 1.0. Further, we assume the weight $k$
of $f$ to be $<0$. For each $\lambda \in \Lambda$, we have 
\begin{eqnarray}\label{first1}
a(\lambda)=\sum_{x} \int_{{\rm C}_x^{*}} f(z) e^{-2\pi i \lambda z} dz,
\end{eqnarray}
where $x$ runs over a system of representatives of $\Q$ modulo $\Z$. Moreover, 
the integrals 
and the series in this expression converge absolutely.}\par
We note that the function $z \mapsto f(z)e^{-2\pi i\lambda z}$ is invariant by
$z \mapsto z+1$. Hence our series does not depend on the choice of the system of 
representatives.
By \eqref{eq2}, we know that 
\begin{equation}\label{1coefficient}
a(\lambda) = \int_{[i, i+1]} f(z)e^{-2\pi i \lambda z} dz.
\end{equation}
Since $f$ is holomorphic on $\goth H$, we can replace the segment of integration $[i, i+1]$
by any other continuous path connecting $i$ to $i+1$ in $\goth H$. \par

We shall construct such a path ${\bf c}_{\rm N}$ for each integer ${\rm N} \geqslant 1$. 
For this purpose, we
consider the ordered sequence of rational numbers in $[0,1]$ whose denominators are less
than or equal to ${\rm N}$: this sequence is called the {\em Farey sequence of 
order ${\rm N}$.}
Two consecutive terms $x' = \frac{a'}{c'}$ and $x = \frac{a}{c}$ of this sequence, 
where the fractions are
in irreducible form with positive denominators, are such that $ac' -a'c =1$. Hence the 
corresponding
Ford circles ${\rm C}_{x'}$ and ${\rm C}_x$ are tangent.
The path ${\bf c}_{\rm N}$ is constructed by starting at $i$ on the first circle ${\rm C}_0$
and moving along it in the clockwise direction till it meets the second circle and then 
moving similarly along the second circle. Continuing this we reach $i+1$ through
an arc of the last Ford circle ${\rm C}_1$ in the clockwise direction. 
(For more details see \cite{rad}, section 116.)\par

When ${\rm N}$ goes to $\infty$, the arc of a given Ford circle ${\rm C}_x$ contained in
${\bf c}_{\rm N}$ increases and its limit is ${\rm C}_x^{*}$ if $x \in ]0,1[$,
the right half of ${\rm C}_0^{*}$ if $x=0$ and the left half of ${\rm C}_1^{*}$ if $x=1$. 
Hence proposition 1 follows from the equality
 $a(\lambda) =  \int_{{\bf c}_{\rm N}} f(z)e^{-2\pi i \lambda z} dz$ 
by passing to the limit, if we are able to prove that the sum 
$\sum_{x \in [0,1[ \cap {\Q}} \int_{{\rm C}_x^{*}} |f(z) e^{-2\pi i \lambda z}| |dz|$ is finite. \par

Let $x$ be a rational number. We can choose 
$\gamma = \big({ a\ b \atop c\ d }\big) \in \SL_2(\Z)$ 
such that $c>0$ and $x = \frac{a}{c}$.
Under the change of variable $z \mapsto \gamma z$, the line $i+ \R$  goes
to ${\rm C}_x^{*}$ and so we have
\begin{eqnarray}\label{changevariable}
\int_{{\rm C}_x^{*}}|f(z)e^{-2\pi i \lambda z}|~|dz| 
&=& \int_{i+\R} |f(\gamma z) e^{-2\pi i \lambda \gamma z}|~|d(\gamma z)| \nonumber \\ 
&=& \int_{i+\R} |j(\gamma, z) f(z) e^{-2\pi i \lambda \gamma z}|~|d(\gamma z)| \\
&=& \int_{i+\R} |cz+d|^{k-2} |f(z) e^{-2\pi i \lambda \gamma z}|~|dz|. \nonumber
\end{eqnarray}
Now the function $f$ is bounded on $i+\R$.  We have
$$
|e^{-2\pi i \lambda \gamma z}|  =  e^{2\pi \lambda {\rm Im}(\gamma z)}
\leqslant ~e^{2\pi \lambda_+ /c^2} \leqslant ~e^{2\pi \lambda_+}
$$ 
for $z \in i+ \R$, where $ \lambda_+ = {\rm max}(0, \lambda)$, and
$$ 
\int_{i+\R} |cz+d|^{k-2}~|dz| = c^{k-2} \int_{i+\R} |z|^{k-2}~|dz|.
$$
Therefore there exists a constant ${\rm A}>0$, independent of $x$,
such that
\begin{eqnarray}\label{second1}
\int_{{\rm C}_x^{*}}|f(z)e^{-2\pi i \lambda z}|~|dz| \leqslant {\rm A} c^{k-2}.
\end{eqnarray}
The sum of these integrals for $x \in [0,1[ \cap {\Q}$ is finite:
indeed, for each $c\geqslant1$, the number of rational numbers in $ [0,1[$ with denominator
$c$ is  $\varphi(c)$ (where $\varphi$ denotes the Euler function), hence is at most $c$;
and since $k<0$, the series $\sum_{c=1}^{\infty} c^{k-1}$ converges. This 
completes the proof.

\medskip
\noindent
{\bf Remark.}$~-$
The previous proof shows that if we assume $k<1$ in place of \hbox{$k<0$}, each integral  
$\int_{{\rm C}_x^{*}} f(z) e^{-2\pi i \lambda z} dz$ still makes sense and is
absolutely convergent.

\subsection*{1.2}{\bf Some special functions} 

\smallskip

To compute the integrals occurring in proposition~1, 
we shall need to study the functions ${\rm L}_\nu$ defined for $\nu$ and $t$ in $\C$ by 
\begin{equation}\label{Bessel}
{\rm L}_\nu (t) = \sum_{n=0}^{\infty}  \frac{t^n}{n! \,\Gamma(n + \nu +1)}\ckomma
\end{equation}
where $\Gamma$ is the Gamma function. Since $1/\Gamma$ is an entire 
function and since the series \eqref{Bessel}
converges normally on any compact set of $\C^2$, the map $(\nu,t) 
\mapsto {\rm L}_\nu(t)$ is holomorphic on 
$\C^2$. The function ${\rm L}_\nu$ is essentially a Bessel function: more precisely, 
the modified Bessel function of the first kind, classically denoted by 
${\rm I}_\nu(t)$, is equal to 
$(\frac{t}{2})^\nu {\rm L}_\nu(\frac{t^2}{4})$. 

\bigskip
\noindent
{\bf Proposition 2.}$-$
{\em For $\nu \in \C$ with ${\rm Re}(\nu) >0$ and $t \in \C$, we have
\begin{eqnarray}\label{Besselintegral}
{\rm L}_\nu(t)  = \frac{1}{2\pi i}\int_{\sigma - i\infty}^{\sigma + i\infty}
w^{-\nu -1} e^{w+\frac{t}{w}} ~dw,
\end{eqnarray}
where $\sigma$ is any positive real number and the complex power $w^{-1-\nu}$ is defined
as in the introduction.}

\medskip
Before proving proposition 2, we shall recall two versions of Hankel's formula
for the inverse of the $\Gamma$ function.

\bigskip
\noindent
{\bf Lemma 1.}$-$
{\em 
$a)$ Let $\sigma$ be a positive real number. We have
\begin{equation}\label{hankel}
\frac{1}{\Gamma(s)} = \frac{1}{2\pi i}\int_{\sigma - i\infty}^{\sigma + i\infty}e^w w^{-s} dw
\end{equation}
for all $s \in \C$ with ${\rm Re}(s) > 1$.\par
$b)$
Let $\sigma$ and $\beta$ be positive real numbers. Let ${\rm H}_{\sigma, \beta}$ denote
the following Hankel contour: it consists of the horizontal half line from
$-\infty -i\beta$ to $\sigma -i\beta$, followed by the vertical segment 
from $\sigma - i\beta$ to $\sigma + i\beta$ and by the horizontal half line
from $\sigma + i\beta$ to $-\infty + i\beta$. We have
\begin{equation}\label{hankel1}
\frac{1}{\Gamma(s)} = \frac{1}{2\pi i}\int_{{\rm H}_{\sigma, \beta}}e^w w^{-s} dw
\end{equation}
for all $s \in \C$.} \par
For $w \in \C \bfminus \R_-$, we have 
\begin{equation}\label{powers}
|w^{-s}| = e^{-{\rm Re}(s\log w)} 
\leqslant e^{-{\rm Re (s)}\log|w| + \pi |{\rm Im}(s)|} 
\leqslant |w|^{-{\rm Re}(s)}e^{\pi |{\rm Im}(s)|}
\end{equation}
and  $|e^w| = e^{{\rm Re}(w)}$. 
It follows that the integral in \eqref{hankel} converges for ${\rm Re}(s) >1$, that the
integral in \eqref{hankel1} converges for all $s\in \C$ and that these two integrals 
are equal when ${\rm Re}(s) >1$. It therefore suffices to prove assertion $b)$.
Since both sides of \eqref{hankel1} are entire functions of $s$, it even suffices
to prove \eqref{hankel1} when ${\rm Re}(s) <0$. In this case, since the right hand side
of \eqref{hankel1} does not depend on $\beta$, we can replace it by its 
limit when $\beta$ goes to $0$. This limit 
is 
$$
\frac{1}{2\pi i} \int_{-\infty}^{0} e^w |w|^{-s} 
\left( e^{s\pi i} - e^{-s\pi i} \right)dw
= \frac{\sin \pi s}{\pi} \int_{0}^{\infty} e^{-w} w^{-s} dw 
= \frac{\sin \pi s}{\pi}\Gamma(1-s).
$$ 
It is equal to $\frac{1}{\Gamma(s)}$ by the reflection formula. This completes the 
proof of the lemma.

\bigskip
We now prove proposition 2. Applying formula \eqref{hankel} for $s = n+ \nu+ 1$, we get 
from \eqref{Bessel} 
$$ 
{\rm L}_\nu(t) = \frac{1}{2\pi i}\sum_{n=0}^{\infty}  
\frac{t^n}{n!}\int_{\sigma - i\infty}^{\sigma + i\infty}e^w w^{-n-\nu-1} dw,
$$
for $\nu\in \C$ such that ${\rm Re}(\nu) >0$ and $t \in \C$.
Proposition 2 follows by interchanging the sum and the integral. In order to 
justify this interchange, 
we shall prove that the integral
$$
\int_{\sigma - i\infty}^{\sigma + i\infty}\sum_{n=0}^{\infty} 
\frac{|t|^n}{n!} |e^w w^{-n-\nu-1}| |dw|
= \int_{\sigma - i\infty}^{\sigma + i\infty} e^{|t|/|w|} |e^w w^{-\nu-1}| |dw|
$$ 
is finite.
The function $w \mapsto e^{|t|/|w|} |e^w|$ is bounded on the line of integration.
By \eqref{powers}, we have $|w^{-\nu-1}| = {\rm O}(|w|^{-{\rm Re}(\nu)-1})$
when $w$ goes to infinity on this line.
Since ${\rm Re}(\nu) > 0$, 
$ \int_{\sigma - i\infty}^{\sigma + i\infty}|w|^{-{\rm Re}(\nu)-1} |dw|$
converges. This completes the proof.

\medskip
\noindent
{\bf Remark.}$~-$
Let ${\rm H}_{\sigma, \beta}$ be the Hankel contour defined in lemma 1, $b)$. We
have 
\begin{eqnarray}\label{BesselHankel}
{\rm L}_\nu(t)  = \frac{1}{2\pi i}\int_{{\rm H}_{\sigma, \beta}}
w^{-\nu -1} e^{w+\frac{t}{w}} ~dw
\end{eqnarray}
for all $\nu \in \C$. The proof is similar to the proof of proposition 2 except that
we use formula \eqref{hankel1} instead of formula \eqref{hankel} and the fact that
$\int_{{\rm H}_{\sigma, \beta}} |w^{-\nu-1}e^{w}|  |dw|$ converges.

\subsection*{1.3}{\bf Computation of integrals} 

\smallskip
In this subsection, we compute the integrals occurring in proposition~1. We adopt
the hypotheses and notations of subsection 1.0, and we moreover assume the weight
$k$ of $f$ to be $<0$. 

\bigskip
\noindent
{\bf Proposition 3.}$-$
{\em Let $\lambda$ be any element in $\Lambda$, $x$ be a rational number
and ${\rm C}_x$ denote the corresponding Ford circle.
Choose $\gamma = \big( {a \ b \atop c \ d}\big) \in \SL_2(\Z)$ with $c>0$ such that
$x = \frac{a}{c}$. Write $j(\gamma,z) = \varepsilon(\gamma) \left(\frac{cz+d}{i}\right)^k$.
The integral $\int_{{\rm C}_x^{*}} f(z)e^{-2\pi i \lambda z} dz$ is equal to
\begin{equation}\label{prop3}
(2\pi)^{2-k}c^{k-2}\sum_{\mu\in\Lambda, \mu<0}a(\mu)|\mu|^{1-k}
\varepsilon(\gamma)e^{-2\pi i \frac{\lambda a+ \mu d}{c}} 
{\rm L}_{1-k}\left(4{\pi}^2\lambda|\mu|c^{-2}\right).
\end{equation}
}\par

Under the change of variable $z \mapsto \gamma z$, the line $i+\R$ goes to
${\rm C}_x^{*}$. The negative orientation of $i+\R$ corresponds 
to the clockwise orientation of ${\rm C}_x^{*}$. Hence we have
\begin{equation}\label{eq4}
\int_{{\rm C}_x^{*}} f(z)e^{-2\pi i \lambda z} dz = \int_{i+\infty}^{i-\infty}
j(\gamma,z)f(z)e^{-2\pi i \lambda \gamma z} d(\gamma z).
\end{equation}

\bigskip
\noindent
{\bf Lemma 2.}$-$
{\em We have 
\begin{equation}\label{eq5}
\int_{i+\infty}^{i-\infty} j(\gamma,z)(f(z) - {\rm P}(z)) 
e^{-2\pi i \lambda \gamma z} d(\gamma z) = 0.
\end{equation}
}
The function $f- {\rm P}$ is bounded on the half plane 
$\overline{\goth H}_1 = \{z \in \C|{\rm Im}(z) \geqslant 1\}$.
The same is true of the function $z \mapsto e^{-2\pi i \lambda \gamma z}$ since 
$$
|e^{-2\pi i \lambda \gamma z}| = e^{ 2\pi \lambda {\rm Im}(z)/|cz+d|^2} 
\leqslant e^{2 \pi \lambda_+/c^2}
$$
for $z \in \overline{\goth H}_1$, where $ \lambda_+ = {\rm max}(0, \lambda)$.
We have $|j(\gamma, z)| |\frac{d(\gamma z)}{dz}| = |cz+d|^{k-2}$ and $k$ 
is negative, hence the left hand side of \eqref{eq5} remains unchanged when we replace
$\int_{i+\infty}^{i-\infty}$ by $\int_{it+\infty}^{it-\infty}$ for any
$t\geqslant1$, and its absolute value is bounded by a constant times
$$
\int_{it+ \R}|cz+d|^{k-2}|dz| = c^{k-2}t^{k-1}\int_{i+ \R} |z|^{k-2} |dz| .
$$
Lemma 2 follows by letting $t$ go to $\infty$. 

\bigskip
From \eqref{eq4} and \eqref{eq5}, we deduce
\begin{eqnarray} \label{sanoli}
&& \int_{{\rm C}_x^{*}} f(z)e^{-2\pi i \lambda z} dz =  \int_{i+\infty}^{i-\infty}
j(\gamma,z){\rm P}(z)e^{-2\pi i \lambda \gamma z} d(\gamma z) \\
&=& \sum_{\mu \in \Lambda, \mu<0}a(\mu) \varepsilon(\gamma)
\int_{i-\infty}^{i+\infty}\big(\frac{cz+d}{i}\big)^{k-2}
e^{2\pi i  (\mu z - \lambda \frac{az+b}{cz+d})} dz \nonumber
\end{eqnarray}
By making the change of variable $w = 2\pi i \mu (z+\frac{d}{c})$, we see that the
last expression is equal to
\begin{eqnarray*}
\sum_{\mu \in \Lambda, \mu<0}a(\mu) \varepsilon(\gamma)
e^{-2\pi i \frac{\lambda a + \mu d}{c}} \big(\frac{c}{2\pi|\mu|}\big)^{k-2}
\int_{2\pi |\mu|-i\infty}^{2\pi |\mu| +
 i\infty}w^{k-2}e^{w +\frac{4\pi^2\lambda |\mu|}{c^2w}}\frac{dw}{2\pi i|\mu|}
\end{eqnarray*}
and by proposition 2, the last integral is equal to $
|\mu|^{-1}{\rm L}_{1-k}(4\pi^2\lambda |\mu|c^{-2})$.
This proves proposition 3.

\bigskip
\noindent
{\bf Remark 1.}$~-$
Let $\lambda, \mu \in \Lambda$. The complex number $\varepsilon(\gamma)
e^{-2\pi i \frac{\lambda a + \mu d}{c}}$ 
has absolute value 1
and depends only on $\lambda, \mu$ and $x$ and not on the choice of the 
matrix $\gamma = \big({ a\ b \atop c\ d }\big)$ in $\SL_2(\Z)$ such that $c>0$ and 
$x=\frac{a}{c}$. 
Moreover, it remains unchanged if we translate $x$ by an integer. These statements
are proved by checking that if we replace $\gamma$ by 
$\gamma \big({ 1\ n \atop 0\ 1}\big)$ or by
$\big({ 1\ n \atop 0\ 1}\big)\gamma$, $\varepsilon(\gamma)$ gets multiplied
by $u^n$ and $e^{-2\pi i \frac{\lambda a + \mu d}{c}}$ by $u^{-n}$, where
$u$ is the constant value taken by $j(\big({ 1\ 1 \atop 0\ 1}\big), z)$.
It will be convenient to adopt the notation $\alpha_x(\lambda, \mu)$ for the
complex number $\varepsilon(\gamma)e^{-2\pi i \frac{\lambda a + \mu d}{c}}$.

\medskip
\noindent
{\bf Remark 2.}$~-$
As we have seen in the remark of subsection 1.1, the integrals 
$\int_{{\rm C}_x^{*}} f(z)e^{-2\pi i \lambda z} dz $
still converge if we assume $k<1$ in place of $k<0$. Proposition 
3 remains valid in this context
with the same proof.

\subsection*{1.4}{\bf Rademacher's formula for the coefficients of $f$} 

\smallskip
Hypotheses and notations are as in 1.0 and we moreover assume the 
weight $k$ of $f$ to be $<0$. 
For any pair $(\lambda, \mu)$ of elements
of $\Lambda$ and any integer $c \geqslant 1$, let
${\rm A}_c(\lambda, \mu)$ denote the sum of the complex numbers $\alpha_x(\lambda, \mu)$
defined in the remark 1 of 1.3, where $x = \frac{a}{c}$ and $a$ runs over the 
integers in $[0, c[$ coprime to $c$. Since this sum has exactly $\varphi(c)$ terms and since 
$|\alpha_x(\lambda, \mu)|=1$ for each of them, we have 
$|{\rm A}_c(\lambda, \mu)| \leqslant \varphi(c) \leqslant c$.

\bigskip
By combining proposition 1 and proposition 3, we get Rademacher's formula
expressing all coefficients of the expansion of $f$ at $\infty$
in terms of those of its polar part.

\bigskip
\noindent
{\bf Theorem 4.}$-$
{\em Let $\lambda \in \Lambda$. With the previous notations, we have  
\begin{equation}\label{rademacher}
a(\lambda) = (2\pi)^{2-k} \sum_{c=1}^{\infty}  c^{k-2}
\sum_{\mu \in  \Lambda, \mu <0~} a(\mu) |\mu|^{1-k} {\rm A}_c(\lambda, \mu)
~{\rm L}_{1-k}\left(4{\pi}^2\lambda|\mu|c^{-2}\right).
\end{equation}
}

\bigskip
\noindent
{\bf Remark 1.}$-$
Theorem 4 applied to $\eta^{-1}$, where $\eta$ is the Dedekind eta 
function yields Rademacher's
formula for the number of partitions $p(n)$ of an integer $n$. Note that in this case
$k$ is equal to $- \frac{1}{2}$ and ${\rm L}_{\frac{3}{2}}(t)$ to 
$\frac{1}{\sqrt{\pi}}\frac{d}{dt}\left(\frac{{\rm sh}(2\sqrt{t})}{\sqrt{t}}\right)$. 

\medskip
\noindent
{\bf Remark 2.}$-$
The automorphy factor $j(\big({ 0\ -1 \atop 1\ \ 0 }\big), z)$ can be written as
$\varepsilon\left(\frac{z}{i}\right)^k$ for some $\varepsilon \in \C$,
and we have $j(\big({ -1\ \ 0 \atop \!\!  \ \ 0 \  -1 }\big),z) =1$. By applying
the cocycle relation~\eqref{cocycle} to $\gamma= \gamma' = 
\big({ 0\ -1 \atop 1\ \ 0 }\big)$ and $z=i$,
we get $\varepsilon^2 =1$. Hence $\varepsilon$ is equal to $1$ or $-1$. 
We note that ${\rm A}_1(\lambda, \mu) = \varepsilon$ for any pair 
$(\lambda, \mu) \in \Lambda^2$,
as one sees by computing $\alpha_0(\lambda, \mu)$ using the matrix
$\gamma = \big({ 0\ -1 \atop 1\ \ 0 }\big)$.

\subsection*{1.5}{\bf Proof of theorem 2 when $\Gamma = \SL_2(\Z)$} 

\smallskip
Assumptions and notations are as in subsection 1.0 and we moreover assume the weight $k$ 
of $f$ to be $<0$. There exists at least one $\mu <0$ in $\Lambda$ such that $a(\mu) \ne 0$,
otherwise $f$ would be equal to $0$ by theorem 4. Let $\mu_0$ be the one with 
largest absolute
value. Let $\delta = |\mu_0|$ and
$\delta'$ denote the second largest of the numbers $|\mu|c^{-2}$, where $\mu$ is a negative
element in $\Lambda$ such that $a(\mu) \ne 0$ and $c$ is a positive integer.

\bigskip
\noindent
{\bf Proposition 4.}$-$
{\em For each $\lambda \geqslant 0$ in $\Lambda$, we have
\begin{eqnarray}\label{prop4}
&& \big| a(\lambda) - (2\pi)^{2-k} a(\mu_0) \delta^{1-k} \varepsilon
~{\rm L}_{1-k}\left(4{\pi}^2\lambda \delta \right) \big| \\
&& \phantom{mmmmmmm} 
\leqslant (2\pi)^{2-k} \zeta(1-k) {\rm M}(f) \delta^{1-k} 
~{\rm L}_{1-k}\left(4{\pi}^2\lambda \delta' \right),
\end{eqnarray}
where $\zeta$ denotes the Riemann zeta function, ${\rm M}(f)$ is equal to 
$\sum_{ \mu \in \Lambda, \mu <0} |a(\mu)|$
and $\varepsilon = \pm 1$ is defined in remark 2 of subsection 1.4.
}\par

We apply theorem 4. By remark 2 of subsection 1.4, the term
in the double sum \eqref{rademacher} corresponding to $c=1$ and $\mu=\mu_0$ is
$(2\pi)^{2-k} a(\mu_0) \delta^{1-k} \varepsilon
~{\rm L}_{1-k}\left(4{\pi}^2\lambda \delta \right)$.
We bound the absolute value 
of all the terms corresponding to the other pairs $(c, \mu)$ such that $a(\mu) \ne 0$ 
by using the following facts: \par
$a)$ By formula \eqref{Bessel}, the function ${\rm L}_{1-k}$ is non-negative and increasing
on $\R_{+}$. Hence 
$0 \leqslant {\rm L}_{1-k}\left(4{\pi}^2\lambda |\mu| c^{-2} \right)
 \leqslant {\rm L}_{1-k}\left(4{\pi}^2\lambda \delta' \right)$.\par
$b)$ We have $|\mu|^{1-k} \leqslant \delta^{1-k}$, $|A_c(\lambda, \mu)|\leqslant c$ and 
$\sum_{c=1}^{\infty} c^{k-1} = \zeta(1-k)$.

\bigskip
\noindent
{\bf Corollary 1.}$-$
{\em The coefficient $a(\lambda)$ is equivalent to 
\begin{eqnarray}\label{cor1}
(2\pi)^{2-k} a(\mu_0) \delta^{1-k} \varepsilon
~{\rm L}_{1-k}\left(4{\pi}^2\lambda \delta \right)
\end{eqnarray}
when $\lambda \in \Lambda$ goes to $+\infty$.
}\par

Since $\varepsilon =\pm 1$ and $a(\mu_0) \ne 0$, the corollary follows 
from proposition 4 if we prove that
${\rm L}_{1-k}\left(4{\pi}^2\lambda \delta' \right)/{\rm L}_{1-k}\left(4{\pi}^2\lambda \delta \right) $ 
goes to $0$ when $\lambda \to +\infty$. This can be deduced from the 
classical fact in the theory of Bessel functions that 
${\rm L}_{1-k}(t)$ is equivalent to $\frac{1}{\sqrt{4\pi}} t^{\frac{2k-3}{4}} e^{2\sqrt{t}}$
when $t \to +\infty$, but it is also a consequence of the following more elementary lemma,
applied to the function ${\rm L}_{1-k}$.

\bigskip
\noindent
{\bf Lemma 3.}$-$
{\em Let ${\rm L} \ne 0$ be an entire function defined by a formal power series with 
non-negative coefficients $\sum_{n=0}^{\infty} b_n t^n$. Let $\alpha'$ and $\alpha$
be two positive real numbers such that $\alpha' < \alpha$. The function
$t \mapsto \frac{{\rm L}(t \alpha')}{{\rm L}(t \alpha)}$ is decreasing on $\R_+^*$.
Its limit when $t \to +\infty$ is $0$ if ${\rm L}$ is not a polynomial.}\par

Let $t'$ and $t$ be two positive real numbers such that $t' \leqslant t$. We have
\begin{eqnarray*}
{\rm L}(t\alpha){\rm L}(t'\alpha') - {\rm L}(t\alpha'){\rm L}(t'\alpha) \hskip -.5cm
&&= \sum_{n,m}b_nb_m\left(t^n{\alpha}^n t'^m{\alpha'}^m 
- t^n{\alpha'}^n {t'}^m{\alpha}^m \right)  \\
&&= \sum_{n,m}b_nb_mt^n{t'}^m\left({\alpha}^n{\alpha'}^m - {\alpha'}^n{\alpha}^m \right)\\
&&= \sum_{n >m} b_nb_m\left(t^n{t'}^m - t^m {t'}^n\right)\left({\alpha}^n{\alpha'}^m 
- {\alpha'}^n{\alpha}^m \right)
\geqslant 0,
\end{eqnarray*}
hence $\frac{{\rm L}(t' \alpha')}{{\rm L}(t'\alpha)} 
\geqslant \frac{{\rm L}(t \alpha')}{{\rm L}(t \alpha)}$.
This proves that the function
$t \to \frac{{\rm L}(t \alpha')}{{\rm L}(t \alpha)}$ is decreasing on~$\R_+^*$. 
Let $\ell$ denote its limit when $t \to +\infty$. If ${\rm L}$ is not 
a polynomial, then for each 
${\rm N} \in \N$,
${\rm L}(t)$ is equivalent to 
${\rm L}_{\rm N}(t) = \sum_{n={\rm N}}^{\infty} b_n t^n$ when $t \to +\infty$.
But ${\rm L}_{\rm N}(t \alpha') \leqslant 
\left(\frac{\alpha'}{\alpha}\right)^{\rm N}{\rm L}_{\rm N}(t \alpha)$
for all $t \in \R_{+}$,
hence $\ell \leqslant \left(\frac{\alpha'}{\alpha}\right)^{\rm N}$. We thus have $\ell =0$.

\bigskip
\noindent
{\bf Corollary 2.}$-$
{\em There exists $\lambda_0 \in \R$ such that $a(\lambda) \ne 0$ 
for all $\lambda \geqslant \lambda_0$
in $\Lambda$.
}

\bigskip
Such a real number $\lambda_0$ can be determined effectively: indeed by proposition 4 
and lemma~3, it suffices
to choose $\lambda_0$ in $\R_+$ such that
\begin{eqnarray}\label{fourth}
\frac{{\rm L}_{1-k}\left(4{\pi}^2\lambda_0 \delta' \right)}{{\rm L}_{1-k}
\left(4{\pi}^2\lambda_0 \delta \right)}
< \zeta(1-k)^{-1}\frac{|a(\mu_0)|}{\sum_{\mu \in \Lambda, \mu <0} |a(\mu)|} \, \cdot
\end{eqnarray}

\section*{{\S 2.} A variant of the circle method}

In this section, we shall present a variant of Rademacher's circle method based on 
the theory of Poincar\'e series. Assumptions and notations are the same as 
in subsection 1.0. 
Moreover, we suppose that the weight $k$ of $f$ is $<0$. Recall that the automorphy
factors of $f$ are denoted by $j(\gamma, z)$ for \hbox{$\gamma \in \SL_2(\Z)$}, that $u$
denotes the constant value of $j(\big({ 1 \ 1 \atop 0  \  1 }\big),z)$ and
$\sum_{\lambda \in \Lambda} a(\lambda) e^{2\pi i\lambda z}$ the expansion of $f$ at 
$\infty$, where $\Lambda = \left\{\lambda \in \C | e^{2\pi i \lambda} =u \right\}$.

\subsection*{2.1}{\bf Principle of the method} 

\smallskip
In subsection 2.2, we shall construct for each $\lambda$ in 
$\Lambda$ a Poincar\'e  series 
${\rm E}_{\lambda}$ with the following properties: \par
$a)$ It is a modular function of weight $2-k$ for $\SL_2(\Z)$, 
holomorphic on $\goth H$, with automorphy factors 
$j(\gamma, z)^{-1}(cz+d)^2$ for $\gamma = \big({ a\ b \atop c\ d }\big) \in \SL_2(\Z)$.
We then have ${\rm E}_{\lambda}(z+1) = u^{-1} {\rm E}_{\lambda}(z)$.\par
$b)$ The function ${\rm E}_{\lambda}(z) - e^{- 2\pi i \lambda z}$ goes to $0$ when
${\rm Im} (z)$ goes to $+\infty$.\par
\noindent
It follows from $a)$ and $b)$ that   ${\rm E}_{\lambda}(z) - e^{- 2\pi i \lambda z}$
has an expansion at $\infty$ of the form 
\begin{eqnarray}\label{fifth}
{\rm E}_{\lambda}(z) - e^{- 2\pi i \lambda z} = \sum_{\mu \in \Lambda, \mu <0} e(\lambda, \mu)e^{-2\pi i \mu z}.
\end{eqnarray}
The coefficients $e(\lambda,\mu)$ will be computed in subsection 2.3.

\medskip

The differential form $f(z){\rm E}_{\lambda}(z)dz$ is then holomorphic on $\goth H$ and
invariant under the action of $\SL_2(\Z)$. It defines by passing to the quotient a
 holomorphic 
differential form on the Riemann surface $\SL_2(\Z) \backslash \goth H$. The
modular invariant identifies this Riemann surface with ${\bf P}_1(\C)\!\bfminus\!\{\infty\}$. 
But any holomorphic differential form on ${\bf P}_1(\C) \!\bfminus \!\{\infty\}$
has residue $0$ at $\infty$. Therefore we have
$$
\int_{[i, i+1]} f(z) {\rm E}_{\lambda}(z)dz = 0.
$$
We use the expansions of $f$ and ${\rm E}_{\lambda}$ at $\infty$ to compute this integral
and we get 
\begin{equation}\label{poincare}
a(\lambda) + \sum_{\mu \in \Lambda, \mu <0} a(\mu) e(\lambda, \mu) = 0.
\end{equation}
Replacing the coefficients $e(\lambda,\mu)$ by their expressions obtained in subsection~2.3, 
we shall obtain the same formulas as in theorem 4 for $\lambda$ in $\Lambda$.

\subsection*{2.2}{\bf The Poincar\'e series ${\rm E}_{\lambda}$} 

\smallskip
Let ${\rm U}$ denote the stabilizer of $\infty$ in $\SL_2(\Z)$: it is generated by the matrices 
$\big({ -1\ \ 0 \atop \ \ 0\ -1 }\big)$  and $\big({ 1\ 1 \atop 0\ 1 }\big)$. We have 
$j(\big({ -1\ \ 0 \atop \ \ 0\ -1 }\big), z) = 1$ and $j( \big({ 1\ 1 \atop 0\ 1 }\big), z) = u$. \par
Let $\lambda \in \Lambda$. Since $e^{2\pi i \lambda} =u$, the differential form
$j(\gamma, z) e^{-2\pi i \lambda \gamma z} d(\gamma z)$ on $\goth H$ depends only on the coset
${\rm U}\gamma$ in ${\rm U} \backslash \SL_2(\Z)$. We define a function ${\rm E}_{\lambda}$
on $\goth H$ by 
\begin{eqnarray}\label{seven}
{\rm E}_{\lambda}(z)dz = \sum_{\gamma \in {\rm U} \backslash \SL_2(\Z)} j(\gamma, z) 
e^{-2\pi i \lambda \gamma z} d(\gamma z)\cdot
\end{eqnarray}
The function ${\rm E}_{\lambda}$ is the sum of $e^{-2\pi i \lambda z}$ and of the series 
\begin{eqnarray}\label{partial}
\sum_{\gamma \in {\rm U} \backslash \left(\SL_2(\Z) \bfminus {\rm U}\right)}
 j(\gamma, z) e^{-2\pi i \lambda \gamma z} \frac{d(\gamma z)}{dz}.
\end{eqnarray}

\medskip
\noindent
{\bf Lemma 4.}$-$
{\em The series \eqref{partial} converges normally on any half plane ${\goth H}_r$ 
with $r>0$ and its sum goes to $0$ when
${\rm Im}(z)$ goes to $+\infty$.}\par

Let $\gamma = \big({ a \ b  \atop c\ d }\big) \in \SL_2(\Z) \bfminus {\rm U}$. 
Then $c \ne 0$
and we have $\big|j(\gamma, z) \frac{d(\gamma z)}{dz} \big|= |cz+d|^{k-2}$ and
$|e^{-2\pi i \lambda\gamma z}| = e^{2\pi \lambda {\rm Im}(z)/|cz+d|^2} \leqslant 
e^{2 \pi \lambda_+ /{\rm Im}(z)}
\leqslant e^{2 \pi \lambda_+ /r}$ on ${\goth H}_r$. \par

The map $\big({ a \ b  \atop c\ d }\big) \mapsto {\rm sign}(c)(c,d)$
defines by passing to the quotient a bijection from ${\rm U} \backslash
\left(\SL_2(\Z) \bfminus {\rm U}\right)$
onto the set of pairs of coprime integers $(c,d)$ such that $c \geqslant 1$.
The sum over these pairs of $|cz+d|^{k-2}$ is finite because
$k < 0$. This proves that the series \eqref{partial} converges normally on ${\goth H}_r$. 
Since each
of its terms goes to $0$ when ${\rm Im}(z)$ goes to $+\infty$, the same is true of
its sum.

\bigskip
It follows from lemma 4 that ${\rm E}_{\lambda}$ is a holomorphic function on 
$\goth H$. It follows
from formula \eqref{seven} and from the cocycle relation \eqref{cocycle} that 
the holomorphic differential form 
$\omega_{\lambda} = {\rm E}_{\lambda}(z) dz$ on $\goth H$ satisfies the functional equation
$\gamma^{*}(\omega_{\lambda}) = j(\gamma, z)^{-1} \omega_{\lambda}$ for $\gamma$ in 
$\SL_2(\Z)$.
Equivalently, we have ${\rm E}_{\lambda}(\gamma z) 
= j(\gamma, z)^{-1}(cz+d)^2 {\rm E}_{\lambda}(z)$
for $\gamma = \big({ a\ b \atop c\ d }\big)$ in $\SL_2(\Z)$ and $z \in \goth H$. By
lemma 4, ${\rm E}_{\lambda}(z) - e^{- 2\pi i \lambda z}$ goes to $0$ when ${\rm Im}(z)$
goes to $\infty$. Furthermore it gets multiplied by $u^{-1}$ when $z$ is replaced 
by $z+1$. Hence its expansion at $\infty$ is of the form 
$\sum_{\mu \in \Lambda, \mu <0} e(\lambda, \mu)e^{-2\pi i \mu z}$.

\subsection*{2.3}{\bf Computation of the coefficients $e(\lambda, \mu)$}

\smallskip
For all $\mu <0$ in $\Lambda$, we have by lemma 4
\begin{eqnarray}\label{coefficientP}
e(\lambda, \mu) &=&  \int_{[i, i+1]} e^{2\pi i \mu z} ({\rm E}_{\lambda}(z) - e^{-2\pi i \lambda z})dz\\
&=& \sum_{\gamma \in {\rm U} \backslash \left(\SL_2(\Z) \bfminus {\rm U}\right)}  \int_{[i, i+1]}
j(\gamma, z) e^{2\pi i (\mu z -\lambda \gamma z)}d(\gamma z).\nonumber
\end{eqnarray}
Let $\gamma$ be an element of $\SL_2(\Z) \bfminus {\rm U}$ and let $\gamma_0$ denote the matrix
$\big({ 1\ 1 \atop 0 \ 1 }\big)$. Then 
${\rm U} \cap \gamma{\rm U}\gamma^{-1} = \left\{\big({ 1\ 0 \atop 0 \ 1 }\big), 
\big({ -1\  \ 0 \atop \ 0 \ -1 }\big)\right\}$.
Therefore the cosets ${\rm U}\gamma\gamma_0^n$, for $n \in \Z$, are all distinct. Their union
is the double coset ${\rm U}\gamma{\rm U}$. Moreover, since $e^{2\pi i \mu} = u$
and $j(\gamma\gamma_0^n, z) = j(\gamma, z+n) u^n$, we have
\begin{eqnarray*}
&&\int_{[i, i+1]} j(\gamma\gamma_0^n, z) e^{2\pi i (\mu z -\lambda \gamma\gamma_0^n z)}d(\gamma\gamma_0^n z)\\
&=& \int_{[i, i+1]} 
j(\gamma, z+n) e^{2\pi i (\mu (z+n) -\lambda \gamma (z+n))}d(\gamma (z+n))\\
&=&\int_{[i+n, i+n+1]}
j(\gamma, z) e^{2\pi i (\mu z -\lambda \gamma z)}d(\gamma z).
\end{eqnarray*}
Hence by grouping in \eqref{coefficientP} the terms by double cosets, we get
$$
e(\lambda, \mu) = \sum_{\gamma \in {\rm U} \backslash \left(\SL_2(\Z) \bfminus {\rm U}\right)/{\rm U}}  
\int_{ i - \infty}^{i+ \infty} j(\gamma, z) e^{2\pi i(\mu z -\lambda \gamma z)}d(\gamma z).
$$
These integrals are exactly those that we already met in formula \eqref{sanoli}.
Moreover, the map $\gamma \mapsto \gamma \infty$ defines by passing to the quotient a bijection 
from $\left(\SL_2(\Z) \bfminus {\rm U}\right)/ {\rm U}$ to $\Q$ and hence from
${\rm U} \backslash \left(\SL_2(\Z) \bfminus {\rm U}\right)/{\rm U}$ to $\Q/\Z$.
The same computation as in subsection 1.3 and 1.4 shows that  
$$
e(\lambda, \mu) = - \sum_{c=1}^{\infty} (2\pi)^{2-k}c^{k-2}|\mu|^{1-k}{\rm A}_c(\lambda, \mu) 
{\rm L}_{1-k}\left(4 \pi^2 \lambda |\mu|c^{-2}\right),
$$
the notations being those of subsection 1.4. This shows formula \eqref{poincare} is
the same as theorem 4.

\section*{{\S 3.} Extension of the circle method to non-negative weight}

\subsection*{3.0}{\bf Assumptions and notations} 

\smallskip

Assumptions and notations in this section are those of subsection 1.0 but we now
restrict ourselves to the case where {\em the weight $k$ of the modular function $f$
is non-negative}. We recall that the automorphy
factors of $f$ are denoted by $j(\gamma, z)$ for $\gamma \in \SL_2(\Z)$, that $u$
denotes the constant value of $j(\big({ 1\ 1 \atop 0  \ 1 }\big),z)$, 
$\sum_{\lambda \in \Lambda} a(\lambda) e^{2\pi i\lambda z}$ the expansion of $f$ at 
$\infty$, where $\Lambda = \left\{\lambda \in \C | e^{2\pi i \lambda} =u \right\}$,
and $\rm P$ denotes its polar part.

\subsection*{3.1}{\bf The circle method revisited} 

\smallskip

Let $\lambda \in \Lambda$. By \eqref{eq2}, we know that
\begin{equation}\label{Lcoefficient}
a(\lambda) = \int_{[i, i+1]} f(z)e^{-2\pi i \lambda z} dz.
\end{equation}
Since $f$ is holomorphic on $\goth H$, we can replace the segment of integration $[i, i+1]$
by any other continuous path connecting $i$ to $i+1$ in $\goth H$, in particular
by the paths ${\bf c}_{\rm N}$ defined in 1.1. But the procedure of letting ${\rm N}$
go to $\infty$ is no more justified. \par
Hence from now on, we fix an integer ${\rm N} \geqslant 1$. Let $\Q_{\rm N}$ denote the set of
rational numbers with denominator $\leqslant {\rm N}$. It is a closed discrete subset of $\R$.
Let $x \in \Q_{\rm N}$. Let $x'$ denote the predecessor and $x''$ the successor of $x$.
The Ford circle ${\rm C}_{x}$ is tangent to ${\rm C}_{x'}$ at a point ${\rm R}_{x, {\rm N}}$
and to ${\rm C}_{x''}$ at a point ${\rm S}_{x, {\rm N}}$. Let  ${\rm C}_{x, {\rm N}}$
denote the arc of ${\rm C}_{x}$ joining ${\rm R}_{x, {\rm N}}$ to ${\rm S}_{x, {\rm N}}$
in the clockwise direction. The path ${\bf c}_{\rm N}$ is composed of the right half of
${\rm C}_{0, {\rm N}}$, of the arcs ${\rm C}_{x, {\rm N}}$ for $x \in ]0,1[$ and of the
left half of ${\rm C}_{1, {\rm N}}$. Hence we have 
\begin{equation}\label{circlemethod}
a(\lambda) = \sum_{x}\int_{{\rm C}_{x, {\rm N}}} f(z)e^{-2\pi i \lambda z} dz,
\end{equation}
where $x$ runs over a system of representatives of ${\Q}_{\rm N}$
modulo $\Z$.

\subsection*{3.2}{\bf Parametrization of the paths ${\rm C}_{x, {\rm N}}$} 

\smallskip
Let $x$ be an element of ${\Q}_{\rm N}$. Choose $\gamma = \big({ a\ b \atop c \ d }\big)$ 
in $\SL_2(\Z)$ with $c >0$ such that $x= \frac{a}{c}$.
Let $n$ and $m$ denote the largest integers in $\Z$ such that $d+nc \leqslant {\rm N}$
and $-d+mc \leqslant {\rm N}$ respectively. We set
$$
a' = b+ na, \quad c' = d+ nc,  \quad a'' = -b+ ma, \quad c'' = -d+ mc.
$$
We have $1 \leqslant c' \leqslant {\rm N}$, $1\leqslant c''
 \leqslant {\rm N}$, $c+ c' \geqslant {\rm N} +1$,
$c + c'' \geqslant {\rm N} + 1$. We hence have $c^2 + {c'}^2 
\geqslant \frac{(c+c')^2}{2} \geqslant \frac{({\rm N} + 1)^2}{2}$
and similarly $c^2 + {c''}^2 \geqslant \frac{({\rm N} + 1)^2}{2}$. With these notations:

\bigskip
\noindent
{\bf Proposition 5.}$-$
{\em $a)$ The predecessor of $x$ in ${\Q}_{\rm N}$ is $x'=\frac{a'}{c'}$ and its successor is  
$x''=\frac{a''}{c''}$. \par
$b)$ We have ${\rm R}_{x, {\rm N}}= \frac{ac + a'c' + i}{c^2 + {c'}^2} = \frac{a'+ ia}{c' + ic}$
and ${\rm S}_{x, {\rm N}}= \frac{ac + a''c'' + i}{c^2 + {c''}^2} = \frac{a'' -  ia}{c'' - ic}$.\par
$c)$ The imaginary parts of ${\rm R}_{x, {\rm N}}$ and ${\rm S}_{x, {\rm N}}$ lie in 
the interval $[\frac{1}{2{\rm N}^2}, \frac{2}{({\rm N} + 1)^2}]$.\par
$d)$ The oriented path ${\rm C}_{x, {\rm N}}$ is the image under the map $z \mapsto \gamma z$
of the oriented segment $[n+i, -m+i]$ of the line $i +\R$.}\par

The point $x'=\frac{a'}{c'}$ belongs to ${\Q}_{\rm N}$
and is smaller than $x$ since $ac'- a'c = ad - bc =1$. Let $\frac{p}{q}$ with $q > 0$ 
and ${\rm gcd}(p,q)=1$ be a rational number such that $x' < \frac{p}{q} < x$. We have $aq -cp \geqslant 1$
and $c'p - a'q \geqslant 1$, hence $q = c'(aq -cp) + c(c'p -a'q) \geqslant c' + c > {\rm N}$.
Therefore $\frac{p}{q}$ does not belong to ${\Q}_{\rm N}$. This
proves that $x'$ is the predecessor of $x$ in ${\Q}_{\rm N}$. The proof that $x''$ is its successor 
is similar.\par

Since the radius of the circle ${\rm C}_{x'}$ is $\frac{1}{2{c'}^2}$ and the radius of 
${\rm C}_{x}$ is $\frac{1}{2{c}^2}$, the point ${\rm R}_{x, {\rm N}}$ is the barycenter
of the centers $\frac{a'}{c'}+ \frac{i}{2{c'}^2}$ and $\frac{a}{c} + \frac{i}{2{c}^2}$
of these circles with coefficients $\frac{1}{2{c}^2}$ and $\frac{1}{2{c'}^2}$, or
equivalently ${c'}^2$ and $c^2$. We hence have 
${\rm R}_{x, {\rm N}}= \frac{ac + a'c' + i}{c^2 + {c'}^2} = \frac{a'+ ia}{c' + ic}$.
Since $\frac{({\rm N} +1)^2}{2} \leqslant c^2 + {c'}^2 \leqslant 2{\rm N}^2$, 
${\rm Im}({\rm R}_{x, {\rm N}}) = \frac{1}{c^2 + {c'}^2}$
belongs to $[\frac{1}{2{\rm N}^2}, \frac{2}{({\rm N} + 1)^2}]$.
This proves $b)$ and $c)$ for ${\rm R}_{x, {\rm N}}$. The proofs for ${\rm S}_{x, {\rm N}}$
are similar. \par

The map $z\mapsto \gamma z$ defines a diffeomorphism from $i +\R$ to ${\rm C}_x \bfminus \{x\}$,
the negative orientation of $i + \R$ corresponding to the clockwise orientation of ${\rm C}_x$.
This map sends $n+i$ to ${\rm R}_{x, {\rm N}}$ and $-m +i$ to ${\rm S}_{x, {\rm N}}$ by $b)$.
This proves $d)$.

\bigskip
\noindent
{\bf Corollary.}$-$
{\em 
$a)$ Both points ${\rm R}_{x, {\rm N}}$  and ${\rm S}_{x, {\rm N}}$ are at a distance
$\leqslant \frac{\sqrt{2}}{c{\rm N}}$ from $x$.\par
$b)$ The image of the chord $[{\rm R}_{x, {\rm N}}, {\rm S}_{x, {\rm N}}]$ of ${\rm C}_x$ by $\gamma^{-1}$ 
is an arc of a circle of radius $\leqslant \frac{{\rm N}^2}{c^2}$.
}\par
We have ${\rm R}_{x, {\rm N}} - x = \frac{a'+ ia}{c' + ic} - \frac{a}{c} = \frac{-1}{c(c' + ic)}$.
The distance between ${\rm R}_{x, {\rm N}}$ and $x$ is therefore $\frac{1}{c\sqrt{c^2 + {c'}^2}}$.
Since $c^2 + {c'}^2 \geqslant \frac{{\rm N}^2}{2}$, this distance is $\leqslant \frac{\sqrt{2}}{c{\rm N}}$.
This proves~$a)$ for ${\rm R}_{x, {\rm N}}$. The proof for ${\rm S}_{x, {\rm N}}$
is similar.

Let ${\rm D}$ denote the line passing through the points ${\rm R}_{x, {\rm N}}$ and ${\rm S}_{x, {\rm N}}$
and $\delta$ denote its distance to $x$. Since ${\rm R}_{x, {\rm N}}$
lies on the left half and ${\rm S}_{x, {\rm N}}$ on the right
half of the circle ${\rm C}_x$, the point of $\rm D$ closest to $x$ belongs to 
the chord $[{\rm R}_{x, {\rm N}}, {\rm S}_{x, {\rm N}}]$. The
ordinate of this point is $\geqslant \frac{1}{2{\rm N}^2}$ by proposition 5, $c)$, 
hence $\delta \geqslant \frac{1}{2{\rm N}^2}$.
The map $z \mapsto \gamma^{-1}z = \frac{dz-b}{-cz +a}= -\frac{d}{c} - \frac{1}{c^2\left(z - \frac{a}{c}\right)}$ 
is obtained by composing the translation $z \mapsto z - \frac{a}{c}$, 
the transformation $z \mapsto -\frac{1}{c^2 z}$ and the translation $z \mapsto z-\frac{d}{c}$.
The first one maps ${\rm D}$ to a line ${\rm D'}$ whose distance to the origin is $\delta$, the second one 
maps ${\rm D'}$ into a circle $\rm C'$ of radius $\frac{1}{2c^2\delta}$ and the third one maps  
$\rm C'$ onto a circle of the same radius. Since $\frac{1}{2c^2 \delta} \leqslant \frac{{\rm N}^2}{c^2}$,
assertion $b)$ follows.

\subsection*{3.3}{\bf First error term: comparing $f$ to its polar part} 

\smallskip
We keep the notations of subsection 3.2. It follows from proposition 5, $d)$ that we have
\begin{equation}\label{segment}
\int_{{\rm C}_{x, {\rm N}}} f(z)e^{-2\pi i \lambda z} dz
= \int_{n+i}^{-m+i} j(\gamma ,z)f(z)e^{-2\pi i \lambda \gamma z} d(\gamma z)
\end{equation}
for all $\lambda \in \Lambda$. 

\medskip
We recall that $\rm P$ denotes the polar part $\sum_{\mu \in \Lambda, \mu <0} 
a(\mu) e^{2\pi i \mu z}$
of $f$ at $\infty$. The function $|f-{\rm P}|$ is bounded on 
$\overline{\goth H}_1 = \{ z \in \C| {\rm Im}(z) \geqslant 1\}$.
We denote by $\| f -{\rm P}\|_{\overline{\goth H}_1}$ its supremum on
$\overline{\goth H}_1$.

\bigskip
\noindent
{\bf Proposition 6.}$-$
{\em
For each $\lambda \in \Lambda$, $|\int_{n+ i}^{-m+i}j(\gamma,z)(f(z)- 
{\rm P}(z))e^{-2\pi i \lambda \gamma z} d(\gamma z)|$
is bounded by 
\begin{equation}\label{firsterror}
\pi e^{4 \pi \lambda_+/({\rm N}+1)^2} \| f -{\rm P}\|_{\overline{\goth H}_1}  \times 
\begin{cases} 
(2{\rm N}^2)^{k-1} c^{-k} & \text{if $k \geqslant 2$,} \\
2{\rm N}^{\frac{3k}{2} -1}c^{-\frac{k}{2} -1} 
& \text{if $0 \leqslant k\leqslant 2$.}  
\end{cases}
\end{equation}
}\par
The integral ${\rm I} =\int_{n+i}^{-m+i}j(\gamma,z)(f(z)- {\rm P}(z))
e^{-2\pi i \lambda \gamma z} d(\gamma z)$
is unchanged if we replace the path of integration $[n+i, -m+i]$ by any other
continuous path from $n+i$ to $-m+i$, for example by $\gamma^{-1}([{\rm R}_{x, {\rm N}}, 
{\rm S}_{x, {\rm N}}])$.
Since $[{\rm R}_{x, {\rm N}}, {\rm S}_{x, {\rm N}}]$ is contained in the closed disc bounded
by ${\rm C}_x$, the path $\gamma^{-1}([{\rm R}_{x, {\rm N}}, {\rm S}_{x, {\rm N}}])$ is 
contained in the
closed half plane $\overline{\goth H}_1$ bounded by $i + \R$. Hence $|f - {\rm P}|$
is bounded by $\| f -{\rm P}\|_{\overline{\goth H}_1}$ on this path. By proposition 5, $c)$,
$|e^{-2\pi i \lambda \gamma z}| = e^{2 \pi \lambda {\rm Im}(\gamma z)}$ is bounded by 
$e^{4 \pi \lambda_+/({\rm N}+1)^2}$ on this path. We thus have
$$
|{\rm I}| \leqslant e^{4 \pi \lambda_+/({\rm N}+1)^2} \| f -{\rm P}\|_{\overline{\goth H}_1}
\int_{\gamma^{-1}([{\rm R}_{x, {\rm N}}, {\rm S}_{x, {\rm N}}])}
|cz+ d|^{k-2}|dz|.
$$
We shall now distinguish two cases: \par
$a)$ We have $k \geqslant 2$. We have
$|cz+d| = |-c(\gamma z)+a|^{-1} \leqslant (c\,{\rm Im}(\gamma z))^{-1} 
\leqslant 2{\rm N}^2c^{-1}$
by proposition 5, $c)$. Furthermore the length of the path 
$\gamma^{-1}([{\rm R}_{x, {\rm N}}, {\rm S}_{x, {\rm N}}])$ 
is at most $2\pi \frac{{\rm N}^2}{c^2}$ by corollary to proposition 5. Hence we have 
$$
\int_{\gamma^{-1}([{\rm R}_{x, {\rm N}}, {\rm S}_{x, {\rm N}}])}
|cz+ d|^{k-2} |dz| \leqslant \pi(2{\rm N}^2)^{k-1} c^{-k}.
$$
\indent
$b)$ We have $0 \leqslant k \leqslant 2$. By H\"older's inequality,
$\int_{\gamma^{-1}([{\rm R}_{x, {\rm N}}, {\rm S}_{x, {\rm N}}])}
|cz+ d|^{k-2} |dz|$ is bounded by
$$
\left(\int_{\gamma^{-1}([{\rm R}_{x, {\rm N}}, {\rm S}_{x, {\rm N}}])}
|cz+ d|^{-2} |dz|\right)^{1-\frac{k}{2}}
 \left(\int_{\gamma^{-1}([{\rm R}_{x, {\rm N}}, {\rm S}_{x, {\rm N}}])}|dz|\right)^{\frac{k}{2}}.
$$ 
Since $|cz+ d|^{-2} |dz| = |d (\gamma z)|$, the first of these two integrals 
is the length of the segment $[{\rm R}_{x, {\rm N}}, {\rm S}_{x, {\rm N}}]$
and the second the length of the arc of circle
$\gamma^{-1}([{\rm R}_{x, {\rm N}}, {\rm S}_{x, {\rm N}}])$.
The first one is bounded by $\frac{2 \sqrt{2}}{c {\rm N}}$ and the second 
one by $2 \pi \frac{{\rm N}^2}{c^2}$
by corollary to proposition 5. Hence we have  
$$
\int_{\gamma^{-1}([{\rm R}_{x, {\rm N}}, {\rm S}_{x, {\rm N}}])}
|cz+ d|^{k-2} |dz|
\leqslant \Big(\frac{2 \sqrt{2}}{c {\rm N}}\Big)^{1-\frac{k}{2}} 
\Big(2 \pi \frac{{\rm N}^2}{c^2}\Big)^{\frac{k}{2}}
\leqslant 2\pi 
{\rm N}^{\frac{3k}{2} -1}c^{-\frac{k}{2} -1}.
$$
This completes the proof.\par

\subsection*{3.4}{\bf Second error term: comparing a segment to a Hankel contour} 

\smallskip
We keep the notations of subsection 3.2. Moreover, we write 
$j(\gamma ,z) = \varepsilon(\gamma)\left(\frac{cz+d}{i}\right)^k$. We then have 
for each $\lambda \in \Lambda$,
\begin{eqnarray}\label{transformation}
&&\int_{n+ i}^{-m+i} j(\gamma ,z){\rm P}(z)e^{-2\pi i \lambda \gamma z}
 d(\gamma z)  \nonumber \\ 
&=& \sum_{\mu \in \Lambda, \mu <0} a(\mu) \varepsilon(\gamma)
\int_{-m+i}^{n+i}\big(\frac{cz+d}{i}\big)^{k-2}
  e^{2\pi i \left(\mu z -\lambda \frac{az+b}{cz+d}\right)} dz\\
&=&  \sum_{\mu \in \Lambda, \mu <0} a(\mu) \varepsilon(\gamma)
e^{-2\pi i \frac{\lambda a + \mu d}{c}}
\int_{c - ic'}^{c + ic''}v^{k-2} e^{\frac{2\pi}{c}\left(|\mu| v 
+ \frac{\lambda}{v} \right)} \frac{dv}{ic}\, ,\nonumber
\end{eqnarray}
the last equality being obtained by the change of variable $v = \frac{cz+d}{i}\, 
\cdot$ In the
next proposition ${\rm L}_{1-k}$ denotes the function defined in subsection 1.2.

\bigskip
\noindent
{\bf Proposition 7.}$-$
{\em Assume $\lambda \geqslant 0$. The integral 
$\int_{c - ic'}^{c + ic''}v^{k-2} e^{\frac{2\pi}{c}\left(|\mu| v 
+ \frac{\lambda}{v} \right)} \frac{dv}{ic}$ differs from
$$
(2\pi)^{2-k}c^{k-2}|\mu|^{1-k}{\rm L}_{1-k}\left(4{\pi}^2\lambda|\mu|c^{-2}\right)
$$
by at most $e^{4 \pi \lambda/ ({\rm N} + 1)^2} {\rm N}^{k-2} g_{k}(\mu)$, where
\begin{equation}\label{bound}
g_k(\mu) = 
\begin{cases} 
2^{k/2}\left(e^{2\pi |\mu|}\big( 1 + \frac{1}{2\pi |\mu|}\big) 
+ \frac{\Gamma(k-1)}{(2\pi |\mu|)^{k-1}}\right) 
& \text{if $k \geqslant 2$,} \\
e^{2\pi |\mu|}\big(4 + \frac{1}{\pi |\mu|}\big) & \text{if $0 \leqslant k\leqslant 2$.}  
\end{cases}
\end{equation}
}\par

Let $\rm H$ denote the Hankel contour ${\rm H}_{c, {\rm N}}$ defined in lemma 1, $b)$. By the
change of variable $w = \frac{2 \pi |\mu| v}{c}$ and the remark of subsection 1.2, we get
\begin{eqnarray*}
\int_{\rm H}v^{k-2} e^{\frac{2\pi}{c}\left(|\mu| v + \frac{\lambda}{v} \right)} \frac{dv}{ic}
&=& (2\pi)^{2-k}c^{k-2}|\mu|^{1-k} \int_{\frac{2\pi |\mu|}{c}{\rm H}}
w^{k-2} e^{\big(w + \frac{4 \pi^2 \lambda |\mu|}{c^2w}\big)} \frac{dw}{2 \pi i}\\
&=& (2\pi)^{2-k}c^{k-2}|\mu|^{1-k}{\rm L}_{1-k}\left(4{\pi}^2\lambda|\mu|c^{-2}\right).
\end{eqnarray*}
In order to prove the proposition, it therefore suffices
to prove that the integral of $|v^{k-2} e^{\frac{2\pi}{c}\left(|\mu| v 
+ \frac{\lambda}{v} \right)} \frac{dv}{c}|$ on the complement ${\rm H}'$ of 
$[c - i c', c + ic'']$ in $\rm H$ is bounded by 
$e^{4 \pi \lambda/ ({\rm N} + 1)^2} {\rm N}^{k-2} g_{k}(\mu)$.\par

We shall first prove that we have 
${\rm Re}(\frac{1}{v}) = \frac{s}{s^2 + t^2} \leqslant \frac{2c}{({\rm N} +1)^2}$
for $v = s+ it$ on ${\rm H}'$.
When $s = c$ and $t \in [c'', {\rm N}]$, our assertion follows from the
inequality $\frac{c}{c^2 + t^2} \leqslant \frac{c}{c^2 + {c''}^2} 
\leqslant \frac{2c}{({\rm N} +1)^2}$.
When $s = c$ and $t \in [-{\rm N}, -c']$, we argue similarly. When $s \in [0,c]$
and $t = \pm {\rm N}$, our assertion follows from the
inequality $\frac{s}{s^2 + {\rm N}^2} \leqslant \frac{c}{c^2 + {\rm N}^2} 
\leqslant \frac{c}{c^2 + {c'}^2} 
\leqslant \frac{2c}{({\rm N} +1)^2}$. Finally our assertion is obvious when $s \leqslant 0$.
Hence we have $|e^{2\pi \frac{\lambda}{cv}}| = e^{2\pi \frac{\lambda}{c} {\rm Re}(\frac{1}{v})}
\leqslant e^{4 \pi \lambda/ ({\rm N} + 1)^2}$ for any $v$ on ${\rm H}'$. Therefore
$\int_{{\rm H}'} |v^{k-2} e^{\frac{2\pi}{c}\left(|\mu| v 
+ \frac{\lambda}{v} \right)} \frac{dv}{c}| 
\leqslant e^{4 \pi \lambda/ ({\rm N} + 1)^2}{\rm J}$, where
${\rm J} =\int_{{\rm H}'} |v|^{k-2} e^{\frac{2\pi}{c}|\mu| {\rm Re}(v)} \frac{|dv|}{c}$. \par

We shall now bound the integral $\rm J$. The segment $[c+ ic'', c+ i{\rm N}]$
has length at most $c$ since we have $c + c'' > {\rm N}$. On this segment, we
have $|v| \leqslant {\rm N}\sqrt{2}$
and \hbox{$|v| \geqslant \sqrt{c^2 + {c''}^2} \geqslant \frac{{\rm N}}{\sqrt{2}}$}. Hence
the contribution of this segment to the integral~$\rm J$ is at most 
$e^{2\pi |\mu|}({\rm N}\sqrt{2})^{k-2}$
when $k \geqslant 2$ and $e^{2\pi |\mu|}\big(\frac{{\rm N}}{\sqrt{2}}\big)^{k-2}$ when 
$0 \leqslant k \leqslant 2$.
The same holds for the segment $[c - i{\rm N}, c -ic']$.\par

For $v = s+ i {\rm N}$ with $s$ in $] -\infty, c]$, we have $|v| 
\geqslant {\rm N}$, $|v| \leqslant {\rm N} \sqrt{2}$
if $s \in [0,c]$ and $|v| \leqslant {\rm sup}({\rm N}\sqrt{2}, |s| \sqrt{2})$ if $s \leqslant 0$.
Hence the contribution of the half line $] -\infty + i {\rm N}, c + i {\rm N}]$ to the
integral $\rm J$ is bounded by
\begin{eqnarray*}
&& \int_{-\infty}^{c} ({\rm N}\sqrt{2})^{k-2} e^{\frac{2\pi}{c}|\mu| s} \frac{ds}{c}
+ \int_{-\infty}^{0} (|s|\sqrt{2})^{k-2} e^{\frac{2\pi}{c}|\mu| s} \frac{ds}{c}\\
&& \phantom {mm} = ({\rm N}\sqrt{2})^{k-2} \frac{e^{2\pi |\mu|}}{2\pi |\mu|}
+ \frac{\Gamma(k-1)}{(2\pi |\mu|)^{k-1}}(c\sqrt{2})^{k-2} 
\end{eqnarray*}
when $k \geqslant 2$ and by
$$
\int_{-\infty}^{c} {\rm N}^{k-2} e^{\frac{2\pi}{c}|\mu| s} \frac{ds}{c}
= {\rm N}^{k-2} \frac{e^{2\pi |\mu|}}{2\pi |\mu|}
$$
when $0 \leqslant k \leqslant 2$. The same holds for the other half line 
$]- \infty - i{\rm N}, c - i {\rm N}]$.
This completes the proof.

\subsection*{3.5}{\bf An approximate value of the coefficients $a(\lambda)$ for
 $\lambda \geqslant 0$ in $\Lambda$} 

\smallskip
Assumptions and notations are those of subsection 3.0. In order to state the main 
result of this subsection,
it will be convenient to define two constants 
which depend only on the modular form $f$ and neither on $\lambda$ nor on $\rm N$.
The first one is ${\rm M}_1(f) = \| f -{\rm P}\|_{\overline{\goth H}_1}$, where
$\rm P$ is the polar part of $f$ (see subsection~3.3). The second one is 
${\rm M}_2(f) = \sum_{ \mu \in \Lambda, \mu <0} |a(\mu)| g_k(\mu)$, where the real numbers
$g_k(\mu)$ are defined in proposition 7.

\medskip
\bigskip
\noindent
{\bf Proposition 8.}$-$
{\em For each $\lambda \geqslant 0$ in $\Lambda$ and each integer ${\rm N} \geqslant 1$, 
we have 
\begin{eqnarray*}
&& \Big |a(\lambda) - (2\pi)^{2-k} \sum_{c=1}^{\rm N}  c^{k-2}
\sum_{\mu \in  \Lambda, \mu <0~} a(\mu) |\mu|^{1-k} {\rm A}_c(\lambda, \mu)
~{\rm L}_{1-k}\left(4{\pi}^2\lambda|\mu|c^{-2}\right) \Big| \phantom{mmm}  \\
&\leqslant &  
e^{4 \pi \lambda/ ({\rm N} +1)^2} 
\Big( \pi {\rm M}_1(f) \times 
\begin{cases} 
(2{\rm N^2})^{k -1} \zeta(k-1) & \text{if $k > 2$} \\
2{\rm N}^2 (1 + \log {\rm N}) & \text{if $k = 2$} \\
\frac{4}{2-k}{\rm N}^{k}  & \text{if $0 \leqslant k < 2$}  
\end{cases} ~+~ {\rm M}_2(f) {\rm N}^k  \Big),
\end{eqnarray*}
where the numbers ${\rm A}_{c}(\lambda, \mu)$ are defined in subsection 1.4.}\par
In the expression \eqref{circlemethod} of $a(\lambda)$, we replace the integrals 
by those given in the left hand side of \eqref{segment}.
We then replace $f$ by its polar part and get an error term \eqref{firsterror}
that we have to sum up when $x$ runs over a system of representatives of ${\Q}_{\rm N}$
modulo $\Z$, for example the rational numbers $\frac{a}{c}$ 
where $1 \leqslant c \leqslant {\rm N}$,
$0 \leqslant a <c$ and ${\rm gcd}(a,c)=1$. We note that to 
each $c \geqslant 1$ correspond at most
$c$ values of $a$. Each integral is then evaluated
using \eqref{transformation} and proposition 7.
We thus get its approximate value with an error term that we have to sum up 
over the same system of representatives of ${\Q}_{\rm N}$
modulo $\Z$. By following this procedure we get the following bound for the left hand side
in proposition 8:
\begin{eqnarray*}
&& \pi e^{4 \pi \lambda/({\rm N}+1)^2} \| f -{\rm P}\|_{\overline{\goth H}_1}  \times 
\begin{cases} 
(2{\rm N}^2)^{k-1} \sum_{c=1}^{\rm N}c^{1-k} & \text{if $k \geqslant 2$} \\
2{\rm N}^{\frac{3k}{2} -1}\sum_{c=1}^{\rm N}
c^{-\frac{k}{2}} & \text{if $0 \leqslant k\leqslant 2$}  
\end{cases} \\
&& + \phantom{m} e^{4 \pi \lambda/ ({\rm N} + 1)^2} {\rm N}^{k-2} \sum_{\mu \in \Lambda, \mu < 0}
\sum_{c=1}^{\rm N}c \,|a(\mu)|\, g_k(\mu).
\end{eqnarray*}
We can slightly simplify this expression by using the following inequalities
\begin{equation*}
\sum_{c=1}^{\rm N} c^{1-k} \leqslant 
\begin{cases} 
\zeta(k-1) & \text{if $k > 2$} \\
1 +\log{\rm N} & \text{if $k = 2$,} \\  
\end{cases} 
\end{equation*}
\begin{equation*}
\sum_{c=1}^{\rm N} c^{-\frac{k}{2}} \leqslant \int_{0}^{\rm N} t^{-\frac{k}{2}} dt
\,=\, \frac{2}{2-k}{\rm N}^{1 -\frac{k}{2}} \quad \hbox{if $0 \leqslant k <2$,}
\end{equation*}
\begin{equation*}
\sum_{c=1}^{\rm N} c \leqslant {\rm N}^{2}.
\end{equation*}
Proposition 8 follows.

\bigskip
\noindent
{\bf Corollary 1.}$-$
{\em For each $\lambda \geqslant 1$ in $\Lambda$, we have 
\begin{eqnarray*}
&& \Big |a(\lambda) - (2\pi)^{2-k} \sum_{c=1}^{[\sqrt{\lambda}]}  c^{k-2}
\sum_{\mu \in  \Lambda, \mu <0~} a(\mu) |\mu|^{1-k} {\rm A}_c(\lambda, \mu)
~{\rm L}_{1-k}\left(4{\pi}^2\lambda|\mu|c^{-2}\right) \Big| \phantom{mmm}  \\
&\leqslant &  
e^{4 \pi} \Big(\pi {\rm M}_1(f) \times 
\begin{cases} 
(2\lambda)^{k -1} \zeta(k-1) & \text{if $k > 2$} \\
\lambda (2 + \log \lambda) & \text{if $k = 2$} \\
\frac{4}{2-k} \lambda^{k/2}  & \text{if $0 \leqslant k < 2$}  
\end{cases} ~+~ {\rm M}_2(f) {\lambda}^{k/2}  \Big).
\end{eqnarray*}
}
We apply proposition 8 with $\rm N = [\sqrt{\lambda}]$. We note that in this case,
we have ${\rm N} + 1 \geqslant \sqrt{\lambda}$, hence 
$e^{4 \pi \lambda/ ({\rm N} +1)^2} \leqslant e^{4 \pi}$.

\medskip
\noindent
{\bf Remark.}$-$
The exponent $k-1$ of $\lambda$ in the error term of corollary 1 for
$k > 2$ is sharp. Indeed, there exist Eisenstein series of weight $k$ with no polar part
for which $|a(\lambda)|$ grows as fast as $\lambda^{k-1}$ when $\lambda \to +\infty$.

\subsection*{3.6}{\bf Proof of theorem 3 when $\Gamma = \SL_2(\Z)$} 

\smallskip
Assumptions and notations are as in subsection 3.0. We recall that the weight $k$ of $f$ is
assumed to be $\geqslant 0$. In this section we furthermore assume that {\em $f$ is not finite 
at $\infty$}. This means that there exists at least one $\mu < 0$ in $\Lambda$ such that
$a(\mu) \ne 0$. Let $\mu_0$ be the one with largest absolute
value. Let $\delta = |\mu_0|$ and
$\delta'$ denote the second largest of the numbers $|\mu|c^{-2}$, where $\mu$ is a negative
element in $\Lambda$ such that $a(\mu) \ne 0$ and $c$ is a positive integer. \par

As we noticed in remark 2 of subsection 1.4, $j(\big({ 0\ -1 \atop 1\ \ 0 }\big), z)$ can be written as
$\varepsilon\left(\frac{z}{i}\right)^k$ for some $\varepsilon$ equal to $1$ or $-1$ and we have
${\rm A}_1(\lambda, \mu) = \varepsilon$ for any pair $(\lambda, \mu) \in \Lambda^2$.

\bigskip
\noindent
{\bf Proposition 9.}$-$
{\em The coefficient $a(\lambda)$ is equivalent to 
$$
(2\pi)^{2-k} a(\mu_0) \delta^{1-k} \varepsilon
~{\rm L}_{1-k}\left(4{\pi}^2\lambda \delta \right)
$$ 
when $\lambda \in \Lambda$ goes to $+\infty$.
}\par
The previous expression is nothing but the term corresponding to $(\mu, c)= (\mu_0, 1)$
in the double sum appearing in the left hand side of corollary 1 to proposition 8.
The error term on the right hand side is bounded by a power of $\lambda$ when
$\lambda$ goes to $\infty$, hence is negligible with respect to 
${\rm L}_{1-k}\left(4{\pi}^2\lambda \delta \right)$. It therefore suffices
to show that the double sum 
\begin{equation}\label{sum}
\sum_{c, \mu } a(\mu)c^{k-2} |\mu|^{1-k} {\rm A}_c(\lambda, \mu)
~{\rm L}_{1-k}\left(4{\pi}^2\lambda|\mu|c^{-2}\right),
\end{equation}
where $c$ runs from $1$ to $[\sqrt{\lambda}]$ and $\mu$ over the negative elements of $\Lambda$
such that $a(\mu) \ne 0$, with $(\mu, c) \ne (\mu_0, 1)$, is also negligible
with respect to ${\rm L}_{1-k}\left(4{\pi}^2\lambda \delta \right)$
when $\lambda$ goes to $\infty$.\par

The function ${\rm L}_{1-k}$ defined in subsection 1.2 by 
${\rm L}_{1-k}(t) = \sum_{n=0}^{\infty}  \frac{t^n}{n! \,\Gamma(n + 2 -k)}$
may have some negative coefficients, but its coefficients for $n \geqslant k-2$
are~\hbox{$\geqslant 0$}. We shall therefore introduce another function
$\widetilde{\rm L}_{1-k}$ defined by
$$
\widetilde{\rm L}_{1-k}(t) = \sum_{n=0}^{\infty}  \frac{t^n}{n! \,|\Gamma(n + 2 -k)|} \,\cdot
$$
For $t \in \R_{+}$, we have $|{\rm L}_{1-k}(t)| \leqslant  \widetilde{\rm L}_{1-k}(t)$.
Moreover, $\widetilde{\rm L}_{1-k}$ is an increasing function on $\R_+$ and 
${\rm L}_{1-k}(t)$ is equivalent to $\widetilde{\rm L}_{1-k}(t)$ when $t$
goes to $\infty$.\par

These properties and the inequalities 
$|{\rm A}_c(\lambda, \mu)| \leqslant c$ and $\sum_{c=1}^{[\sqrt{\lambda}]}c^{k-1} \leqslant {\lambda}^{k/2}$
allow us to bound the absolute value of \eqref{sum} by 
$$
\lambda^{k/2} \sum_{\mu \in \Lambda, \mu < 0} a(\mu)|\mu|^{1-k} 
~\widetilde{\rm L}_{1-k}\left(4{\pi}^2\lambda \delta'\right).
$$
To complete the proof, we need only to show that $\lambda^{k/2}~\widetilde{\rm L}_{1-k}\left(4{\pi}^2\lambda \delta'\right)$, or equivalently $\lambda^{k/2}~{\rm L}_{1-k}\left(4{\pi}^2\lambda \delta'\right)$,
is negligible with respect to ${\rm L}_{1-k}\left(4{\pi}^2\lambda \delta\right)$.
But this follows from the fact that ${\rm L}_{1-k}(t)$ is equivalent to
$\frac{1}{\sqrt{4\pi}} t^{\frac{2k-3}{4}} e^{2\sqrt{t}}$
when $t \to +\infty$.

\bigskip
\noindent
{\bf Corollary.}$-$
{\em There exists $\lambda_0 \in \R_+$ such that $a(\lambda) \ne 0$ for all $\lambda \geqslant \lambda_0$
in $\Lambda$.
}

\bigskip
If needed, such a $\lambda_0$ could be determined effectively by using the explicit bounds
from corollary 1 to proposition 8.

\section*{{\S 4.} The case of an arbitrary finite index subgroup of $\SL_2(\Z)$}

\subsection*{4.0}{\bf Assumptions and notations}

\smallskip
In this section, we shall indicate how to extend the results of the previous sections
when in place of $\SL_2(\Z)$ one considers one of its finite index subgroups~$\Gamma$. So let $f$
be a non-zero modular function for $\Gamma$ of some real weight~$k$, holomorphic
on $\goth H$, and let $j(\gamma, z)$ for $\gamma$ in $\Gamma$ 
denote its automorphy factors. We have 
\begin{eqnarray}\nonumber
f(\gamma z) = j(\gamma,z) f(z)  \quad \hbox{and} \quad |j(\gamma,z)| = |cz+d|^k  
\end{eqnarray}
for  $\gamma = \big({ a\ b \atop c\ d }\big)$ in $\Gamma$  and $z$ in $\goth{H}$.
The functions $j(\gamma,z)$ satisfy the cocycle relation 
$j(\gamma \gamma', z) = j(\gamma, \gamma' z)j(\gamma', z)$ for $\gamma, \gamma'$ in $\Gamma$
and $z$ in $\goth H$.\par

We denote by $h$ the smallest positive integer such that $\big({ 1\ h \atop 0\ 1 }\big) \in \Gamma$.
The function $j(\big({ 1\ h \atop 0\ 1 }\big),z)$ takes a constant value $u$. We have $|u| =1$ and
$f(z+h) =uf(z)$ for $z \in {\goth H}$. We write
$$
f(z)=\sum_{\lambda \in \Lambda} a(\lambda) e^{2 \pi i\lambda z}
$$
the expansion of $f$ at $\infty$,
where $\Lambda$ is the set of $\lambda \in \R$ such that $e^{2\pi i\lambda h} =u$.
For each $\lambda$ in $\Lambda$, we have
\begin{equation}\label{hcoefficient}
a(\lambda) = h^{-1} \int_{[i,i+h]} f(z) e^{-2\pi i \lambda z} dz.
\end{equation}

\subsection*{4.1}{\bf Further notations}

\smallskip
Let $\gamma = \big({ a\ b \atop c\ d }\big)$ be a matrix in $\SL_2(\Z)$ with $c>0$. 
Define $f_{\gamma}$ to be the function $z \mapsto (\frac{cz+d}{i})^{-k}f(\gamma z)$.
It is a modular function of weight $k$ for $\gamma^{-1}\Gamma\gamma$, with respect
to some automorphy factors. We denote by $h_{\gamma}$ the smallest positive
integer such that $f_{\gamma}(z + h_{\gamma}) = u_{\gamma} f_{\gamma}(z)$ for some
$u_{\gamma}$ in $\C$. We have $|u_{\gamma}|=1$. We denote by
 $\sum_{\lambda \in \Lambda_{\gamma}} a_{\gamma}(\lambda)e^{2\pi i \lambda z}$
the expansion of $f_{\gamma}$ at $\infty$, where $\Lambda_{\gamma}$ is the set of $\lambda$
in $\R$ such that $e^{2\pi i \lambda h_{\gamma}} = u_{\gamma}$, and by ${\rm P}_{\gamma}$
its polar part $\sum_{\lambda \in \Lambda_{\gamma}, \lambda <0} a_{\gamma}(\lambda)e^{2\pi i \lambda z}$.\par

We make some elementary remarks on the dependence of all these data on~$\gamma$:\par
$a)$ Let $\gamma' = \big({ a'\ b' \atop c'\ d' }\big)$ be another matrix in $\SL_2(\Z)$ with $c'>0$,
such that $\Gamma\gamma = \Gamma\gamma'$. Since $f$ is a modular function of weight $k$
for $\Gamma$, $f_{\gamma'}$ differs from $f_{\gamma}$ by a multiplicative scalar
of absolute value $1$. Therefore $|f_{\gamma}|$, $|f_{\gamma} - {\rm P}_{\gamma}|$, $h_{\gamma}$, $u_{\gamma}$, $\Lambda_{\gamma}$
and the numbers $|a_{\gamma}(\lambda)|$ for $\lambda \in \Lambda_{\gamma}$ depend only on the
coset $\Gamma\gamma$.\par 

$b)$ In particular, let $\gamma' = \big({ 1 \ mh \atop  0 \ \ 1 \ }\big)\gamma$ for some
$m \in \Z$. We have $f_{\gamma'} = u^m f_{\gamma}$ and therefore $a_{\gamma'}(\lambda) = u^m a_{\gamma}(\lambda)$
for all $\lambda \in \Lambda_{\gamma}$. \par
$c)$ Let $n \in \Z$ and $\gamma' = \gamma \big({ 1\ n \atop  0\ 1 }\big)$. We have
$f_{\gamma'}(z) = f_{\gamma}(z + n)$, therefore we have $h_{\gamma'} = h_{\gamma}$,
$u_{\gamma'} = u_{\gamma}$, $\Lambda_{\gamma'} = \Lambda_{\gamma}$
and $a_{\gamma'}(\lambda) = e^{2\pi i \lambda n}a_{\gamma}(\lambda)$ for
$\lambda \in \Lambda_{\gamma}$.

\subsection*{4.2}{\bf Extension of the circle method to $\Gamma$}

\smallskip
The following proposition generalizes proposition 1 of subsection 1.1. 

\bigskip
\noindent
{\bf Proposition 10.}$-$
{\em Assumptions and notations are those of subsection 4.0. We furthermore assume the weight $k$
of $f$ to be $<0$. For each $\lambda \in \Lambda$, we have 
$$
a(\lambda)= h^{-1} \sum_{x} \int_{{\rm C}_x} f(z) e^{-2\pi i \lambda z} dz,
$$
where $x$ runs over a system of representatives of $\Q$ modulo $h\Z$
and ${\rm C}_x$ is the Ford circle associated to $x$. Moreover, the integrals 
and the series in this expression converge absolutely.}\par

The proof is similar to the proof of proposition 1 with the two following minor changes: \par
$a)$ One starts with formula \eqref{hcoefficient} in place of formula \eqref{1coefficient}.\par

$b)$ For each $x \in \Q$, one chooses a matrix $\gamma=\big({ a\ b \atop c\ d }\big)$ 
in $\SL_2(\Z)$ such that $c>0$ and
$x= \frac{a}{c}$. Formula \eqref{changevariable} gets then replaced by  
$$
\int_{{\rm C}_x} |f(z) e^{-2\pi i \lambda z}| |dz|
= \int_{i+\R} |cz+d|^{k-2}|f_{\gamma}(z) e^{-2\pi i \lambda \gamma z}|~|dz|.
$$
and one notes that the function $f_{\gamma}$ is bounded on $i + \R$
and depends only on the coset $\Gamma\gamma$ by subsection 4.1, assertion $a)$.

\bigskip
Let $\lambda \in \Lambda$ and $x \in \Q$. Let $\gamma = \big({ a\ b \atop c\ d }\big)$ 
be a matrix in $\SL_2(\Z)$ such that $x = \frac{a}{c}$ and $c>0$. As in subsection 1.3,
one proves successively the following equalities
\begin{eqnarray*}
&&\int_{{\rm C}_x} f(z) e^{-2\pi i \lambda z} dz
= \int_{i+ \infty}^{i-\infty} f(\gamma z) e^{-2\pi i \lambda \gamma z} d(\gamma z) \\
&&= \int_{i- \infty}^{i+ \infty} \big(\frac{cz+d}{i}\big)^{k-2}f_{\gamma}(z) e^{-2\pi i \lambda \gamma z} dz
 = \int_{i- \infty}^{i+ \infty} \big(\frac{cz+d}{i}\big)^{k-2}{\rm P}_{\gamma}(z) e^{-2\pi i \lambda \gamma z} dz\\
&& = \sum_{\mu \in \Lambda_{\gamma}, \mu <0} a_{\gamma}(\mu) \int_{i- \infty}^{i+ \infty} \big(\frac{cz+d}{i}\big)^{k-2} e^{2\pi i (\mu z - \lambda \frac{az+b}{cz+d})} dz\\
&& = (2\pi)^{2-k} c^{k-2} \sum_{\mu \in \Lambda_{\gamma}, \mu <0} a_{\gamma}(\mu) 
e^{-2\pi i  \frac{\lambda a + \mu d }{c}} |\mu |^{1-k} {\rm L}_{1-k}\big(4 \pi^2 \lambda |\mu | c^{-2}\big).
\end{eqnarray*}
It follows from assertion $c)$ in subsection 4.1 that the number $a_{\gamma}(\mu) 
e^{-2\pi i  \frac{\lambda a + \mu d }{c}}$ depends only on $x$ and not on the choice of $\gamma$.
We shall therefore denote it by $\beta_x(\lambda,\mu)$. Since $\Lambda_{\gamma}$ too depends only on $x$,
it is appropriate to denote it by $\Lambda_x$. It follows from assertion $b)$ in subsection 4.1 that
$\beta_x(\lambda,\mu)$ and $\Lambda_x$ depend only on $x$ modulo $h\Z$. \par

Replacing the integrals in proposition 10 by the previous expressions, we get the 
following generalization of theorem 4:

\medskip
\noindent
{\bf Theorem 5.}$-$
{\em 
Let $\lambda \in \Lambda$. With the previous notations, we have
$$
a(\lambda) = (2\pi)^{2-k} h^{-1} \sum_x c^{k-2} \sum_{\mu \in \Lambda_{x}, \mu <0}
\beta_{x}(\lambda,\mu) |\mu |^{1-k} {\rm L}_{1-k}\big(4 \pi^2 \lambda |\mu | c^{-2}\big),
$$
where $x$ runs over the rational numbers in $[0,h[$ and $c$ denotes the
denominator of $x$.}\par

\subsection*{4.3}{\bf Proof of theorem 2}

\smallskip
For each $\alpha >0$, let $\Phi_{\alpha}$ denote the set of pairs $(x,\mu)$ with
the following properties:\par

$a)$ $x$ is a rational number in $[0,h[$ and $\mu$ is a negative element in $\Lambda_x$
such that $|\mu| c^{-2} = \alpha$, where $c \geqslant 1$ denotes the denominator of $x$; \par
$b)$ if $\gamma = \big({ a\ b \atop c\ d }\big)$ is a matrix in $\SL_2(\Z)$ such that 
$x = \frac{a}{c}$, then $a_{\gamma}(\mu) \ne 0$.\par
\noindent
We note that condition $b)$ does not
depend on the choice of the matrix $\gamma$: this follows from assertion $c)$ in
subsection 4.1.

\bigskip
Each set $\Lambda_x$ has only finitely many negative elements and $\Lambda_x$
depends only on the orbit $\Gamma x$ of $x$ under $\Gamma$, by assertion $a)$ of subsection 4.1.
Therefore each set $\Phi_{\alpha}$ is finite and 
$\{\alpha \in  \R_+^* | \Phi_{\alpha} \ne \varnothing \}$ is a closed discrete bounded subset of~$\R_+^*$.
It is non-empty, otherwise $f$ would be equal to $0$ by theorem 5; it is in fact infinite. Let
$\delta$ denote its largest and $\delta'$ its second largest element.

\bigskip
\noindent 
{\bf Lemma 5.}$-$
{\em 
There exists $\lambda$ in $\Lambda$ such that the complex number
$$
{\rm B}(\lambda) = \sum_{(x, \mu) \in \Phi_{\delta}} c^{k-2}\, \beta_x(\lambda, \mu)\, |\mu|^{1-k},
$$
where $c$ denotes the denominator of $x$, is different from $0$.}\par

Let $\lambda_1 \in \Lambda$. Then $\Lambda$ is equal to $\lambda_1 +  \Z h^{-1}$.
For $\lambda = \lambda_1 + n h^{-1}$ in $\Lambda$, we have
\begin{equation}\label{coefficientG}
{\rm B}(\lambda) = \sum_{(x, \mu) \in \Phi_{\delta}} c^{k-2}\, \beta_x(\lambda_1, \mu)\, 
|\mu|^{1-k} e^{-2\pi i nx/h}.
\end{equation}
At least one of the coefficients $\beta_x(\lambda_1, \mu)$ appearing in this expression is 
different from $0$, by definition of $\delta$. Furthermore the characters 
$n \mapsto e^{-2\pi i nx/h}$ of $\Z$, for $x \in [0,h[ \cap \Q$, are distinct, hence
linearly independent over $\C$. Therefore the expression \eqref{coefficientG} cannot
vanish for all $n$ in $\Z$.

\bigskip
\noindent
{\bf Remark.}$-$
Let $\rm D$ denote a common denominator of the rational numbers $x$, for $(x,\mu)$ in $\Phi_{\delta}$.
It follows from formula \eqref{coefficientG} that ${\rm B}(\lambda)$ depends only on $\lambda$
modulo ${\rm D}\Z$. 

\bigskip
The sum $\sum_{\mu \in \Lambda_{\gamma}} |a_{\gamma}(\mu)|\, |\mu|^{1-k}$, for 
$\gamma = \big({a \ b \atop c \ d }\big)$ in $\SL_2(\Z)$ with $c >0$, depends only
on the coset $\Gamma\gamma$ by assertion $a)$ of subsection 4.1. Let ${\rm M}_3(f)$
denote the supremum of these sums.
In the same way as we deduced proposition 4 from theorem 4, we deduce the following result
from theorem 5:

\bigskip
\noindent 
{\bf Proposition 11.}$-$
{\em 
For each $\lambda \geqslant 0$ in $\Lambda$, we have
\begin{eqnarray*}
&& \Big| a(\lambda) - (2\pi)^{2-k} h^{-1} {\rm B}(\lambda) {\rm L}_{1-k}\big(4 \pi^2 \lambda \delta \big)\Big| \\
&&\phantom {mmm} \leqslant (2\pi)^{2-k} \zeta(1-k) {\rm M}_3(f) {\rm L}_{1-k}\big(4 \pi^2 \lambda \delta' \big).
\end{eqnarray*}
} 

\bigskip
By lemma 5, there exists $\lambda_0$ in $\Lambda$
such that ${\rm B}(\lambda_0) \ne 0$. By the remark, we have ${\rm B}(\lambda) = {\rm B}(\lambda_0)$
for all $\lambda \in \lambda_0 + {\rm D}\Z$. Since 
${\rm L}_{1-k}\big(4 \pi^2 \lambda \delta' \big)/{\rm L}_{1-k}\big(4 \pi^2 \lambda \delta \big)$
goes to $0$ when $\lambda \to +\infty$, proposition 11 implies that $a(\lambda)$ is equivalent
to $(2\pi)^{2-k} h^{-1} {\rm B}(\lambda) {\rm L}_{1-k}\big(4 \pi^2 \lambda \delta \big)$
when $\lambda \in \lambda_0 + {\rm D}\Z$ goes to $+\infty$, hence does not vanish
for $\lambda \in \lambda_0 + {\rm D}\Z$ large enough. This proves theorem 2. \par

\subsection*{4.4}{\bf Proof of theorem 3}

\smallskip
Assumptions and notations are those of subsections 4.0 and 4.1. We moreover assume the weight
$k$ of $f$ to be $\geqslant 0$. We fix an integer ${\rm N} \geqslant 1$ and denote by ${\Q}_{\rm N}$
the set of rational numbers with denominator $\leqslant {\rm N}$. \par

For each $\lambda \in \Lambda$, we have
$$
a(\lambda)= h^{-1} \sum_{x} \int_{{\rm C}_{x, \rm N}} f(z) e^{-2\pi i \lambda z} dz,
$$
where $x$ runs over a system of representatives of ${\Q}_{\rm N}$ modulo $h\Z$,
and where ${\rm C}_{x, {\rm N}}$ denotes the arc of the Ford circle ${\rm C}_x$
considered in subsection 3.1. The proof of this formula is the same as the proof of 
\eqref{circlemethod} except that we start from formula \eqref{hcoefficient} in place of
\eqref{Lcoefficient}.

\bigskip
Let $x$ be an element in ${\Q}_{\rm N}$. Choose $\gamma = \big({ a\ b \atop c\ d }\big)$ 
in $\SL_2(\Z)$ such that $x = \frac{a}{c}$ and $c>0$. As in subsection 3.2, let
$n$ and $m$ denote the largest integers in $\Z$ such that $d + nc \leqslant {\rm N}$
and $-d+mc \leqslant {\rm N}$ respectively. By following the same steps as in subsection 3.3, 
one successively proves the following results:\par 
$\rm(i)$ We have for each $\lambda \in \Lambda$
\begin{eqnarray*}
\int_{{\rm C}_{x, \rm N}} f(z) e^{-2\pi i \lambda z} dz 
&&= \int_{n+i}^{-m+i} f(\gamma z) e^{-2\pi i \lambda \gamma z} d(\gamma z)\\
&&= \int_{-m+i}^{n+i}  \big(\frac{cz+d}{i}\big)^{k-2}f_{\gamma}(z) 
e^{-2\pi i \lambda \gamma z} dz.
\end{eqnarray*}\par
$\rm(ii)$ For each $\lambda \in \Lambda$, the integral
$\int_{-m+i}^{n+i}  \big(\frac{cz+d}{i}\big)^{k-2}f_{\gamma}(z) e^{-2\pi i \lambda \gamma z} dz$
differs from $\int_{-m+i}^{n+i}  \big(\frac{cz+d}{i}\big)^{k-2}
{\rm P}_{\gamma}(z) e^{-2\pi i \lambda \gamma z} dz$
by at most 
\begin{equation*}
\pi e^{4 \pi \lambda_+/({\rm N}+1)^2} \| f_{\gamma} -{\rm P}_{\gamma}\|_{\overline{\goth H}_1} 
\times 
\begin{cases} 
(2{\rm N}^2)^{k-1} c^{-k} & \text{if $k \geqslant 2$,} \\
2{\rm N}^{\frac{3k}{2} -1}c^{-\frac{k}{2} -1} 
& \text{if $0 \leqslant k\leqslant 2$.}  
\end{cases}
\end{equation*}\par
$\rm(iii)$ For each $\lambda \in \Lambda$, we have
$$
\int_{-m+i}^{n+i}  \big(\frac{cz+d}{i}\big)^{k-2}{\rm P}_{\gamma}(z) 
e^{-2\pi i \lambda \gamma z} dz
= \sum_{\mu \in \Lambda_{\gamma}, \mu <0} \beta_x(\lambda,\mu) 
\int_{c - ic'}^{c + ic''}v^{k-2} e^{\frac{2\pi}{c}\left(|\mu| v + \frac{\lambda}{v} \right)} \frac{dv}{ic},
$$
where $\beta_x(\lambda,\mu) = a_{\gamma}(\mu)e^{-2\pi i  \frac{\lambda a + \mu d }{c}}$
as in subsection 4.2 and where $c' = d+ nc$ and $c'' = -d + mc$ are the denominators of the predecessor
and of the successor of $x$ in ${\Q}_{\rm N}$ respectively. Moreover proposition 7 yields when
$\lambda \geqslant 0$ the approximate value $(2\pi)^{2-k} c^{k-2}|\mu |^{1-k} 
{\rm L}_{1-k}\big(4 \pi^2 \lambda |\mu | c^{-2}\big)$ 
for the last integral, together with an upper bound for the error term. 

\medskip
The numbers $\| f_{\gamma} - {\rm P}_{\gamma} \|_{\overline{\goth H}_1}$, for 
$\gamma= \big( {a \ b \atop c \ d}\big)$ in $\SL_2(\Z)$ with $c>0$, depend only on the
coset $\Gamma\gamma$ by assertion $a)$ of subsection 4.1. Let ${\rm M}_4(f)$
denote the supremum of these numbers. We define in a similar way ${\rm M}_5(f)$
to be the supremum of the numbers $\sum_{\mu \in \Lambda_{x}, \mu <0} |a_{\gamma}(\mu)| g_k(\mu)$,
with $g_k(\mu)$ as in proposition 7. Summing up the previous results, we get the 
following analog of proposition 8, with the same proof.

\medskip
\noindent
{\bf Proposition 12.}$-$
{\em 
For each $\lambda \geqslant 0$ and each integer ${\rm N} \geqslant 1$, we have
\begin{eqnarray*}
&& \Big| a(\lambda) - (2\pi)^{2-k} h^{-1} \sum_x c^{k-2} \sum_{\mu \in \Lambda_{x}, \mu <0}
\beta_{x}(\lambda,\mu) |\mu |^{1-k} {\rm L}_{1-k}\big(4 \pi^2 \lambda |\mu | c^{-2}\big)\Big|\\
&&\, \leqslant 
e^{4 \pi \lambda/ ({\rm N} +1)^2} 
\Big( \pi {\rm M}_4(f) \times 
\begin{cases} 
(2{\rm N^2})^{k -1} \zeta(k-1) & \text{if $k > 2$} \\
2{\rm N}^2 (1 + \log {\rm N}) & \text{if $k = 2$} \\
\frac{4}{2-k}{\rm N}^{k}  & \text{if $0 \leqslant k < 2$}  
\end{cases} ~+~ {\rm M}_5(f) {\rm N}^k  \Big),
\end{eqnarray*}
where $x$ runs over a system of representatives of 
${\Q}_{\rm N}$ modulo $h\Z$ and $c$ denotes the denominator
of $x$.} \par

Define $\delta$, and $\rm B(\lambda)$ for $\lambda \in \Lambda$, as in subsection 4.3. By
taking $\rm N = [\sqrt{\lambda}]$ in proposition 12 one gets, in exactly 
the same way as in subsection 3.6, the following analog of proposition 9.

\bigskip
\noindent
{\bf Proposition 13.}$-$
{\em 
We have
$$
 \Big| a(\lambda) - (2\pi)^{2-k} h^{-1} {\rm B}(\lambda) 
{\rm L}_{1-k}\big(4 \pi^2 \lambda \delta \big)\Big| 
={\rm o}\big( {\rm L}_{1-k}\big(4 \pi^2 \lambda \delta \big)\big),
$$
when $\lambda$ goes to $+\infty$.} \par

By lemma 5, there exists $\lambda_0$ in $\Lambda$ such that $\rm B(\lambda_0) \ne 0$
and by the remark of subsection 4.3, there exists an integer ${\rm D} \geqslant 1$
such that ${\rm B}(\lambda) = {\rm B}(\lambda_0)$ for $\lambda \in \lambda_0 + {\rm D}\Z$.
It then follows from proposition 13 that $a(\lambda) \ne 0$ for 
$\lambda \in \lambda_0 + {\rm D}\Z$
large enough. This proves theorem 3. \par

\section*{{\S 5.} Lacunarity implies holomorphy on $\goth H$}

Let $g$ be a meromorphic function on the upper-half plane $\goth H$.
Suppose that: \par 
$a)$ $g$ is holomorphic on some half-plane ${\goth H}_r = \{z \in \C |{\rm Im}(z) > r\}$;\par
$b)$ there exists a real number $h > 0$ and a complex number $u$ such that $|u|=1$ 
and $g(z+h) = u g(z)$ for $z \in \goth H$;\par
$c)$ $g$ has finite order at $\infty$.

\bigskip
\noindent
{\bf Proposition 14.}$-$
{\em If the expansion of $g$ at $\infty$ is lacunary, then $g$ 
is holomorphic on $\goth H$.}\par

The expansion of $g$ at $\infty$ can be written in the form
$\sum_{n=0}^{\infty} a(n) e^{2 \pi i (\lambda_0 + n)z/h}$ for some $\lambda_0 \in \R$. 
We can replace $g$ by the function $z \mapsto e^{ -2\pi i \lambda_0 z}g(hz)$
and we are then reduced to the case where $g(z+1) = g(z)$ for $z \in \goth H$
and where the expansion of $g$ at $\infty$ is $\sum_{n=0}^{\infty} a(n) e^{2 \pi i n z}$.\par

In that case, $g$ can be written as $z \mapsto {\rm G}(e^{2 \pi i z})$, where ${\rm G}$
is a meromorphic function on the open disc ${\rm D} = \{z \in \C ~|~ |z| < 1\}$,
holomorphic near $0$ and with Taylor expansion $\sum_{n=0}^{\infty} a(n) q^n$
at $0$.\par

Suppose that $\rm G$ is not holomorphic on $\rm D$. Let $r$
denote the smallest absolute value of the poles of $\rm G$ and let 
$\alpha_1, \cdots, \alpha_m$ be the poles of $\rm G$ with absolute value~$r$,
repeated with multiplicity. The Taylor expansion of 
${\rm G} \prod_{i=2}^{m} (q - \alpha_i)$ around $0$ is
still lacunary. Replacing $\rm G$ by this function, we are reduced
to the case where $\rm G$ has only one pole $\alpha$ with absolute value $r$,
of multiplicity $1$. \par

Let $\beta$ denote the residue of $\rm G$ at this pole. Then
${\rm G} - \frac{\beta}{q-\alpha}$ has no pole of absolute value $r$.
Therefore its radius of convergence is $ > r$.
The coefficients of its Taylor series  
$\sum_{n=0}^{\infty} \left(a(n) + \frac{\beta}{{\alpha}^{n+1}}\right) q^n$
are such that $a(n) + \frac{\beta}{{\alpha}^{n+1}} = {\rm o}(r^{-n})$.
Hence $a(n)$ is equivalent to $\frac{-\beta}{{\alpha}^{n+1}}$
when $n \to +\infty$ and is different from $0$ for $n$ large enough.
This contradicts the fact that the series $\sum_{n=0}^{\infty} a(n) q^n$ is lacunary.

\bigskip
\noindent{\bf Acknowledgement:} The first draft of this paper was written while 
the first author was visiting the Institut de Math\'ematiques de Jussieu in Paris. The final
 draft was written while
the second author was visiting the Institute of Mathematical Sciences in Chennai. 
They thank these institutes for providing excellent working conditions and 
IFIM for financial support.

\end{document}